\documentclass[11pt,letterpaper]{article}
\usepackage{fullpage}
\usepackage{amsmath,amsthm,amsfonts,amssymb}

\newtheorem{theorem}{Theorem}[section]
\newtheorem{lemma}[theorem]{Lemma}

\newtheorem{claim}[theorem]{Claim}

\newtheorem{corollary}[theorem]{Corollary}

\newtheorem{definition}[theorem]{Definition}

\newtheorem{fact}[theorem]{Fact}

\newtheorem{remark}{Remark}[section]

\newcommand{\jnote}[1]{}

\newcommand{\E}{{\mathbb E}}
\newcommand{\e}{\varepsilon}

\newcommand{\dist}{\mathsf{dist}}
\newcommand{\diam}{\mathrm{diam}}
\newcommand{\supp}{\mathrm{supp}}
\newcommand{\remove}[1]{}
\newcommand{\1}{\mathbf{1}}

\newcommand{\eps}{\varepsilon}

\newcommand{\vcon}{\mathsf{con}}

\newcommand{\inter}{\mathsf{inter}}

\renewcommand{\P}{\mathcal P}
\renewcommand{\Pr}{\mathbb P}

\newcommand{\comment}[1]{}

\newcommand{\pieceww}[4]{\left\{\begin{array}{@{}cl@{}}{#1}&{#2}\\{#3}&{#4}\end{array}\right.}

\begin{document}

\titlepage

\title{{\bf Metric uniformization and spectral bounds for graphs}}
\author{Jonathan A. Kelner\thanks{Research partially supported by NSF grant CCF-0843915.} \\ MIT \and
James R. Lee\thanks{Research partially supported by NSF grants CCF-0915251 and CCF-0644037, and a Sloan Research Fellowship.  Part of this work was completed during a visit of the author to the Massachusetts Institute of Technology.}
\\ University of Washington \and
Gregory N. Price\thanks{Research partially supported by NSF grant CCF-0843915, an NSF Graduate Research Fellowship, and an Akamai Fellowship.}
\\ MIT\and
Shang-Hua Teng\thanks{Research partially supported by NSF grant CCF-0635102.}
\\ University of Southern California}

\date{}

\maketitle

\begin{abstract}
We present a method for proving upper bounds on the eigenvalues
of the graph Laplacian.
A main step involves choosing an appropriate ``Riemannian'' metric
to uniformize the geometry of the graph.
In many interesting cases, the existence of such a metric
is shown by examining the combinatorics of special types of flows.
This involves proving new inequalities
on the crossing number of graphs.

In particular, we use our method to show that for any positive integer $k$,
the $k^{\rm th}$ smallest eigenvalue of the Laplacian on an $n$-vertex, bounded-degree
planar graph is $O(k/n)$.  This bound is asymptotically tight for every $k$, as it
is easily seen to be achieved for square planar grids.
We also extend this spectral result to graphs with bounded genus,
and graphs which forbid fixed minors.
Previously, such spectral upper bounds were only known for the case $k=2$.
\end{abstract}

\section{Introduction}

Eigenvalues of the Laplacian on graphs and manifolds have been studied
for over forty years in combinatorial optimization and
geometric analysis.
In combinatorial optimization, spectral methods are a class
  of techniques that use the eigenvectors of
  matrices associated with the underlying graphs.
These matrices include the adjacency matrix, the Laplacian, and
  the random-walk matrix of a graph.
One of the earliest applications of spectral methods is to graph
  partitioning, pioneered by Hall \cite{Hall} and
  Donath and Hoffman \cite{DH72,DH73} in the early 1970s.
The use of the graph Laplacian for partitioning was introduced by Fiedler
  \cite{Fi73,Fi75a,Fi75b}, who showed a connection between the
  second-smallest eigenvalue of the Laplacian of a graph and its
  connectivity.
Since their inception, spectral methods have been
  used for solving a wide range of optimization problems, from graph
  coloring \cite{AspvallGilbert,AlonKahale} to image segmentation
  \cite{ShiMalik,TolliverMiller} to web search \cite{Kleinberg,BrinPage}.

\medskip
\noindent {\bf Analysis of the Fiedler value.}
In parallel with the practical development of spectral methods,
progress on the mathematical front has been extremely fruitful,
involving a variety of connections between various graph properties
and corresponding graph spectra.

In 1970, independent from the work of Hall and of Donath and Hoffman,
  Cheeger \cite{Cheeger} proved that the
  isoperimetric number of a continuous manifold can be
  bounded from above by the square root of the smallest non-trivial
  eigenvalue of its Laplacian.
Cheeger's inequality was then extended to graphs by Alon
  \cite{Alon}, Alon and Milman \cite{AlonMilman}, and Sinclair and
  Jerrum \cite{SinclairJerrum}.
They showed that if the Fiedler value of a graph --- the second
  smallest eigenvalue of the Laplacian of the graph --- is small, then
  partitioning the graph according to the values of the vertices in
  the associated eigenvector will produce a cut where the
  ratio of cut edges to the number of vertices in the cut is similarly small.

Spielman and
  Teng \cite{SpielmanTeng} proved a spectral theorem for planar graphs,
  which asserts that the Fiedler value of every bounded-degree planar
  graph with $n$ vertices is $O (1/n)$.
They also showed that the Fiedler value of a finite-element mesh in
  $d$ dimensions with $n$ vertices is $O (n^{-2/d})$.
Kelner \cite{Kelner} then proved that the Fiedler value of a
  bounded-degree graph with $n$ vertices and genus $g$ is $O((g+1)/n)$.
The proofs in \cite{SpielmanTeng,Kelner} critically use the
   inherent geometric  structure of the planar graphs, meshes, and
   graphs with bounded genus.
Recently, Biswal, Lee, and Rao \cite{BLR} developed a new approach
  for studying the Fiedler value; they resolved most of the open problems in
  \cite{SpielmanTeng}.
In particular, they proved that the Fiedler value of a bounded-degree
graph on $n$ vertices without a $K_{h}$ minor is $O ((h^{6}\log
  h)/n)$.
These spectral theorems together with Cheeger's inequality
  on the Fiedler value immediately imply that one can use
  the spectral method to produce a partition as good as the best known
  partitioning methods for planar graphs \cite{LiptonTarjan},
  geometric graphs \cite{MillerTengThurstonVavasis}, graphs with bounded genus
  \cite{GilbertHutchinsonTarjan}, and graphs free of small complete minors
  \cite{AlonSeymourThomas}.

\medskip
\noindent
{\bf Higher eigenvalues and our contribution.}
Although previous work in the graph setting focuses mostly on $k=2$ (the Fiedler
  value of a graph), higher eigenvalues and eigenvectors are used in
many heuristic algorithms
  \cite{AlpertYaoMoreBetter,ChanSchlagZien,ChanGilbertTeng,TolliverMiller}.
\remove{
Spielman and Teng \cite{SpielmanLect} showed that the embedding of a planar
  graph with the eigenvectors associated with the second and the third smallest
  eigenvalues of the graph Laplacian usually gives a nearly planar
  drawing; and
Tolliver \cite{Tolliver} shows experimentally that higher eigenvectors provide
  better multiway image segmentations.}

In this paper, we prove the following theorem on higher graph spectra,
  which concludes a long line of work on upper bounds for the
  eigenvalues of planar graphs.

\begin{theorem}[Planar and bounded-genus graphs]\label{thm:planar}
Let $G$ be a bounded-degree $n$-vertex planar graph.
Then the $k^{\rm th}$ smallest eigenvalue of
  the Laplacian on $G$ is $O(k/n)$.

More generally, if $G$ can be embedded on an orientable surface
of genus $g$, then the $k^{\rm th}$ smallest eigenvalue of the Laplacian is at most
\begin{equation}
\label{eq:genusbound}
O\left((g+1) (\log (g+1))^2 \frac{k}{n}\right).
\end{equation}
\end{theorem}

The asymptotic dependence on $k$ and $n$ is seen to be tight even for
the special case of square planar grids; see Remark \ref{rem:tightness}.
Our spectral theorem provides a mathematical justification
  of the experimental observation that
  when $k$ is small, the $k^{\rm th}$ eigenvalues of the graphs
  arising in many application domains are small as well.
  We hope our result will lead to new progress in the analysis
  of spectral methods.

  We remark that the $(\log(g+1))^2$ factor of \eqref{eq:genusbound}
  comes from a certain geometric decomposability property of genus-$g$
  graphs (see Theorem \ref{thm:LSgenus}) and is most likely non-essential.
  Without this factor, the bound is tight up to a universal constant,
  as shown by the construction of \cite{GilbertHutchinsonTarjan}.

  \medskip

A well-known generalization of graphs which can be drawn on a manifold
of fixed genus involves the notion of a graph minor.
Given finite graphs $H$ and $G$, one says that $H$ is a {\em minor of $G$}
if $H$ can be obtained from $G$ by a sequence of edge contractions
and vertex deletions.  A family $\mathcal F$ of graphs is said to be {\em minor-closed}
if whenever $G \in \mathcal F$ and $H$ is a minor of $G$, then $H \in \mathcal F$ as well.
By the famous graph minor theorem of Robertson and Seymour \cite{RS04},
every such family $\mathcal F$ is characterized by a finite list of forbidden minors.
For instance, by Wagner's theorem \cite{Wagner}, the family of planar graphs
is precisely the family of graphs which do not have $K_{3,3}$ or $K_5$ as a minor.
We prove the following (see the end of Section \ref{sec:graphs}).

\begin{theorem}[Minor-closed families]
\label{thm:minorclosed}
If $\mathcal F$ is any minor-closed family of graphs which does not contain all graphs,
then there is a constant $c_{\mathcal F} > 0$ such that for all $G \in \mathcal F$ with
$n$ vertices and
maximum degree $d_{\max}$, and all $1 \leq k \leq n$,
$$
\lambda_k(G) \leq c_{\mathcal F} \cdot d_{\max} \frac{k}{n}\,.
$$
\end{theorem}

\medskip
\noindent
{\bf The Riemannian setting and conformal uniformization.}
The spectra of the Laplacian on compact Riemannian surfaces of fixed genus
is also well-studied.  Let $M$ be a compact Riemannian manifold of genus $g$,
and let $\lambda_k(M)$ be the $k^{\rm th}$ smallest eigenvalue
of the Laplace operator on $M$.\footnote{In Riemannian geometry,
the convention is to number the eigenvalues starting from $\lambda_0$,
but we use the graph theory convention to make direct comparison easier.}

Hersch \cite{Hersch70} showed that $\lambda_2(M) \leq O(1/{\rm vol}(M))$ for
Riemannian metrics on the 2-sphere, i.e. for the $g=0$ case.
This was extended by Yang and Yau \cite{YangYau} to a bound
of the form $\lambda_2(M) \leq O((g+1)/n)$ for all $g \geq 0$.
Yau asked  whether, for every $g \geq 0$, there was a constant $c_g$ such that
\begin{equation}\label{eq:korevaar}
\lambda_k(M) \leq c_g \frac{k}{{\rm vol}(M)},
\end{equation}
for all $k \geq 1$.
The question was resolved by Korevaar \cite{Korevaar} who proved
that one can take $c_g = O(g+1)$.
As mentioned at the end of the section, we prove that
bounds in the graph setting yield bounds in the setting of surfaces,
and thus
our result also gives a new proof of \eqref{eq:korevaar}
with the slightly worse constant $c_g = O((g+1) (\log (g+1))^2)$.

\medskip

An important point is that the bounds for planar and genus-$g$ graphs---in addition to the work discussed above for Riemannian surfaces---are
proved using some manifestation of conformal uniformization.
In the graph case, this is via the Koebe-Andreev-Thurston circle packing theorem, and
in the manifold case, by the uniformization theorem.
The methods of Hersch, Yang-Yau, and Spielman-Teng start with a representation
of the manifold or graph on the 2-sphere, and then apply an appropriate
M\"obius transformation to obtain a test vector that bounds $\lambda_2$.
There is no similar method known for bounding $\lambda_3$,
and indeed Korevaar's approach to \eqref{eq:korevaar} is significantly
more delicate and uses very strongly the geometry of the standard 2-sphere.

However, the spectra of graphs may be more subtle than the spectra of
  surfaces.
  We know of a reduction in only one direction:  Bounds on graph eigenvalues
  can be used to prove bounds for surfaces; see Section \ref{sec:surfaces}.
For graphs with large diameter, the analysis of graph spectra resembles
  the analysis for surfaces.
For example, Chung \cite{Chung} gave an upper bound of $O (1/D^{2})$ on the
  the Fiedler value, where $D$ is the diameter of the graph.
Grigor'yan and Yau \cite{GrigoryanYau} extended Korevaar's
  analysis to bounded genus graphs that have a strong volume measure ---
  in particular, these graphs have diameter $\Omega (\sqrt{n})$.

Bounded-degree planar graphs (and bounded genus graphs), however,
  may have diameter as small as $O (\log n)$, making
  it impossible to directly apply these diameter-based spectral
  analyses.
  Our work builds on the method of Biswal, Lee, and Rao \cite{BLR},
  which uses multi-commodity flows to define a deformation of the graph geometry.
  Essentially, we try to construct a metric on the graph
  which is ``uniform'' in a metrically defined sense.
  We then show that sufficiently uniform metrics allow us to
  recover eigenvalue bounds.

  To construct metrics with stronger uniformity properties, which
  can be used to capture higher eigenvalues, we study a new flow problem, which we
  define in Section \ref{sec:prelim} and call {\em subset flows};
  this notion may be independently interesting.  As we discuss in the next section,
  these flows arise as dual objects of certain kinds of optimal spreading metrics
  on the graph.  We use techniques from the theory of discrete metric spaces
  to build test vectors from spreading metrics, and we develop new combinatorial
  methods to understand the structure of optimal subset flows.

  \medskip

Our spectral theorem not only provides a discrete analog
  for Korevaar's theorem on higher eigenvalues, but also extends
  the higher-eigenvalue bounds to graphs with a bounded forbidden minor,
  a family that is more combinatorially defined.
Because the Laplacian of a manifold can be approximated by that of a
sufficiently fine mesh graph (see Section \ref{sec:surfaces}), our result also provides a new proof of
  Korevaar's theorem, with a slightly worse constant.


\subsection{Outline of our approach}\label{sec:outline}

For the sake of clarity,
we restrict ourselves for now to a discussion
of the case where $G=(V,E)$ is
a bounded-degree planar graph.  Let $n=|V|$,
and for $1 \leq k \leq n$, let $\lambda_k$ be
the $k^{\rm th}$ smallest eigenvalue of the Laplacian on $G$
(see Section \ref{sec:prelim-laplacian} for a discussion
of graph Laplacians).
We first review the known methods
for bounding $\lambda_2=\lambda_2(G)$.

\medskip
\noindent
{\bf Bounding $\lambda_2$.}
By the variational characterization of eigenvalues,
giving an upper bound on $\lambda_2$ requires
finding a certain kind of mapping of $G$ into the
real line (see Section \ref{sec:prelim}).
Spielman and Teng \cite{SpielmanTeng} obtain
an initial geometric representation
using the Koebe-Andreev-Thurston circle
packing theorem for planar graphs.
Because of the need for finding a test vector
which is orthogonal to the first eigenvector (i.e., the constant function),
one has to post-process this representation
before it will yield a bound on $\lambda_2$.
They use a topological argument to show the
existence of an appropriate M\"obius transformation
which achieves this.  (As we discussed, a similar step was used
by Hersch \cite{Hersch70} in the
manifold setting.)
Even in the arguably simpler setting of manifolds,
no similar method is known for bounding $\lambda_3$,
due to the lack of a rich enough family
of circle-preserving transformations.

\medskip

Our approach begins with the
arguments of Biswal, Lee, and Rao \cite{BLR}.
Instead of finding an external geometric representation,
those authors begin by finding an appropriate
intrinsic deformation of the graph,
expressed via a non-negative vertex-weighting
$\omega : V \to [0,\infty)$, which induces
a corresponding shortest-path metric\footnote{Strictly speaking, this is only a pseudometric
since $\dist_{\omega}(u,v)=0$ is possible for $u \neq v$, but we ignore this distinction
for the sake of the present discussion.} on $G$,
$$
\mathsf{dist}_{\omega}(u,v) = \textrm{length of shortest $u$-$v$ path},
$$
where the length of a path $P$ is given by $\sum_{v \in P} \omega(v)$.
The proper deformation $\omega$ is found via variational methods,
by minimizing the ratio,
\begin{equation}\label{eq:firstomega}
\frac{\sqrt{\sum_{v \in V} \omega(v)^2}}{\sum_{u,v \in V} \dist_{\omega}(u,v)}.
\end{equation}
The heart of the
analysis involves studying the geometry of the minimal solutions,
via their dual formulation in terms of certain kinds of
multi-commodity flows.
Finally, techniques from the theory of metric embeddings
are used to embed the resulting metric space $(V,\mathsf{dist}_{\omega})$ into
the real line, thus recovering an appropriate test vector to bound $\lambda_2$.

\medskip
\noindent
{\bf Controlling $\lambda_k$ for $k \geq 3$.}
In order to bound higher eigenvalues, we need to produce
a system of many linearly independent test vectors.  The first
problem one encounters is that the optimizer of \eqref{eq:firstomega}
might not contain enough information to produce more than a single
vector if the geometry of the $\omega$-deformed graph
is degenerate, e.g. if $V = C \cup C'$
for two large, well-connected pieces $C,C'$ where $C$ and $C'$ are far apart,
but each has small diameter.  (Intuitively,
there are only two degrees of freedom,
the value of the eigenfunction on $C$ and the value on $C'$.)

\medskip
\noindent
{\bf Spreading metrics and padded partitions.}
To combat this, we would like to impose the constraint
that no large set collapses in the metric $\mathsf{dist}_{\omega}$,
i.e.\ that for some $k \geq 1$ and any subset $C \subseteq V$ with $|C| \geq n/k$, the
diameter of $C$ is large.  In order to produce such an $\omega$
by variational techniques, we have to specify this constraint
(or one like it) in a convex way.  We do this using the
{\em spreading metric} constraints which are well-known in mathematical optimization (see, e.g. \cite{Spreading00}).
The spreading constraint on a subset $S \subseteq V$ takes the form,
\begin{equation}\label{eq:secondomega}
\frac{1}{|S|^2} \sum_{u,v \in S} \dist_{\omega}(u,v) \geq \e \sqrt{\sum_{u \in V} \omega(u)^2},
\end{equation}
for some $\e > 0$.

Given such a spreading weight $\omega$ for sets of size $\approx n/k$,
we show in Section \ref{sec:eigenv-spre-weights} how to obtain a bound
on $\lambda_k$ by producing $k$ smooth, disjointly suppported bump functions on $(V,\mathsf{dist}_{\omega})$,
which then act as our $k$ linearly independent test vectors.
The bound depends on the value $\varepsilon$ from \eqref{eq:secondomega},
as well as a certain geometric decomposability property of the space $(V,\dist_{\omega})$.
The bump functions are produced using padded metric partitions (see, e.g. \cite{KLMN05} and \cite{LN05}),
which are known to exist for all planar graphs from the seminal work
of Klein, Plotkin, and Rao \cite{KPR}.

\medskip
\noindent
{\bf The spreading deformation, duality, and subset flows.}
At this point, to upper bound $\lambda_k$, it suffices to find a spreading weight $\omega$
for subsets of size $\approx n/k$, with $\varepsilon$ (from \eqref{eq:secondomega}) as large as possible.
To the end,
in Section \ref{sec:spre-weights-subset-flows}, we write a convex program
whose optimal solution yields a weight $\omega$
with the largest possible value of $\varepsilon$.
The dual program
involves a new kind of multi-commodity flow problem, which we now describe.

Consider a probability distribution $\mu$ on subsets $S \subseteq V$.
For a flow $F$ in $G$ (see Section \ref{sec:prelim} for a review
of multi-commodity flows), we write $F[u,v]$ for the total amount
of flow sent from $u$ to $v$, for any $u,v \in V$.  In this case
a {\em feasible $\mu$-flow} is one which satisfies, for every $u,v \in V$,
$$
F[u,v] = \Pr_{S \sim \mu}[u,v \in S],
$$
where we use the notation $S \sim \mu$ to denote that $S$ is chosen according to the distribution $\mu$.
In the language of demands, every set $S$ places a demand of $\mu(S)$
between every pair $u,v \in S$.   For instance, the classical
all-pairs multi-commodity flow problem would be specified
by choosing $\mu$ which concentrates all its weight on the entire
vertex set $V$.

Given such a $\mu$, the corresponding ``subset flow'' problem
is to find a feasible $\mu$-flow $F$ so that the total $\ell_2$-norm
of the congestion of $F$ at vertices is minimized (see Section \ref{sec:spre-weights-subset-flows}
for a formal definition of the $\ell_2$-congestion).
Finally, by duality, bounding $\lambda_k$ requires
us to prove {\em lower bounds} on the congestion of {\em every possible}
$\mu$-flow with $\mu$ concentrated on sets of size $\approx n/k$.

\medskip
\noindent
{\bf An analysis of optimal subset flows:  New crossing number inequalities.}
In the case of planar graphs $G$, we use a randomized
rounding argument to relate the existence of
a feasible $\mu$-flow in $G$ with small $\ell_2$-congestion to the ability
to draw certain kinds of graphs in the plane without too many edge crossings.
This was done in \cite{BLR}, where the relevant combinatorial problem
involved the number of edge crossings necessary to draw dense
graphs in the plane, a question which was settled by Leighton \cite{Leighton},
and Ajtai, Chv\'atal, Newborn, and Szemer\'edi \cite{acns}.

In the present work, we have to develop new crossing {\em weight} inequalities
for a ``subset drawing'' problem.
Let $H=(U,F)$ be a graph with non-negative edge
weights $W : F \to [0,\infty)$.  Given a drawing of $H$ in the plane, we define
the crossing weight of the drawing as the total weight of all edge crossings,
where two edges $e,e' \in F$ incur weight $W(e) \cdot W(e')$ when they cross.
Write $\mathsf{cr}(H;W)$ for the minimal crossing weight needed to draw $H$ in the plane.
In Section \ref{sec:subset-congestion}, we prove a generalization of the following theorem
(it is stated there in the language of flows), which
forms the technical core of our eigenvalue bound.

\begin{theorem}[Subset crossing theorem]\label{thm:technicalmain}
There exists a constant $C \geq 1$ such that
if $\mu$ is any probability distribution on subsets of $[n]$ with
$\mathbb E_{S \sim \mu} |S|^2 \geq C$, then the following holds.
For $u,v \in [n]$, let $$W(u,v) = \Pr_{S \sim \mu}[u,v \in S]\,.$$
Then we have,
$$
\mathsf{cr}(K_n;W) \gtrsim \frac{1}{n} \left( \mathbb E_{S \sim \mu} |S|^2\right)^{5/2},
$$
where $K_n$ is the complete graph on $\{1,2,\ldots,n\}$.
\end{theorem}

Observe that the theorem is asymptotically tight for all values of $\mathbb E|S|^2$.
It is straightforward that one can draw an $r$-clique
in the plane using only $O(r^4)$ edge crossings.
Thus if we take $\mu$ to be uniform on $k$ disjoint subsets of size $n/k$, then
the crossing weight is on the order of
$k \cdot (1/k)^2 \cdot (n/k)^4 = n^4/k^5$, which matches the
lower bound $\frac{1}{n} (\mathbb E|S|^2)^{5/2} = \frac{1}{n} (n/k)^5$.
The proof involves
some delicate combinatorial and analytic arguments,
and is discussed at the beginning of Section \ref{sec:subset-congestion}.

\medskip
\noindent
{\bf More general families:  Bounded genus and excluded minors.}
Clearly the preceding discussion was specialized to planar graphs.
A similar approach can be taken for graphs of bounded genus (orientable
or non-orientable) using the appropriate generalization of Euler's formula.

To handle general minor-closed families,
we can no longer deal with the notion of drawings, and
we have to work directly with multi-commodity flows in graphs.  To do this,
we use the corresponding ``flow crossing'' theory
developed in \cite{BLR}, with some new twists to handle the regime
where the total amount of flow being sent is very small (this happens
when bounding $\lambda_k$ for large values of $k$, e.g. $k \geq \sqrt{n}$).

\subsection{Preliminaries}
\label{sec:prelim}

We often use the asymptotic notation $A \lesssim B$ to denote $A=O(B)$.
We use $A \asymp B$ to denote the conjunction of $A \lesssim B$ and $A \gtrsim B$.
For a graph $G$, we use $V(G)$ and $E(G)$ to denote the edge and vertex sets of $G$, respectively.
We write $\mathbb R_+ = [0,\infty)$.

\subsubsection{Laplacian spectrum}\label{sec:prelim-laplacian}
Let $G=(V,E)$ be a finite, undirected graph.
We use $u \sim v$ to denote $\{u,v\} \in E$.
We consider the linear
space $\mathbb R^V = \{ f : V \to \mathbb R \}$ and define the
Laplacian $\mathcal L : \mathbb R^V \to \mathbb R^V$ as the symmetric,
positive-definite linear operator given by
$$ (\mathcal L f)(v) = \sum_{u : u \sim v} \left(f(v) - f(u)\right),$$
which in matrix form could be written as $\mathcal L = D - A$ where
$A$ is the adjacency matrix of $G$ and $D$ the diagonal matrix whose
entries are the vertex degrees.
We wish to give upper bounds on the $k^{\rm th}$ eigenvalue of
$\mathcal L$ for each $k$.  To do this we consider the seminorm given by
$$\| f \|_{\mathcal L}^{2}
   = \langle f, \mathcal L f \rangle = \sum_{u \sim v} (f(u) - f(v))^2,$$
and restrict it to $k$-dimensional subspaces $U \subset \mathbb R^V$.
By the spectral theorem, the maximum ratio $\|f\|_{\mathcal L}^{2} /
\|f\|^{2}$ over $U$ is minimized when $U$ is spanned by the $k$
eigenvectors of least eigenvalue, in which case its value is $\lambda_k$.
Therefore if we exhibit a $k$-dimensional subspace $U$ in which
$\|f\|_{\mathcal L}^{2} \leq c$ for all unit vectors $f$, it follows
that $\lambda_k \leq c$.
In particular, this yields the following simple lemma.

For a function $f : V \to \mathbb R$, write $\supp(f) = \{ x \in V : f(x) \neq 0 \}$.
Furthermore, for $S \subseteq V$, let $N(S) = \{ u \in V : u \in S \textrm{ or } u \sim v \textrm{ for some } v \in S \}$.

\begin{lemma}\label{lem:simple}
For any $k \geq 1$, suppose that $f_1, f_2, \ldots, f_{k} \in \mathbb R^V$
is a collection of non-zero vectors such that for all $1 \leq i, j \leq k$ with $i \neq j$,
$\supp(f_i) \cap N(\supp(f_j)) = \emptyset$.  Then,
$$
\lambda_k \leq \max \left\{\frac{\|f_i\|_{\mathcal L}^{2}}{\|f_i\|^{2}}: i \in \{1,2,\ldots,k\}\right\}.
$$
\end{lemma}

\begin{proof}
Since $\supp(\mathcal L f) \subseteq N(f)$, our assumptions immediately imply that $\langle f_i, \mathcal L f_j\rangle = 0$
for $i \neq j$.
\end{proof}


\subsubsection{Flows}\label{sec:prelim-flows}
Let $G=(V,E)$ be a finite, undirected graph, and for every pair $u,v
\in V$, let $\mathcal P_{uv}$ be the set of all paths between $u$ and
$v$ in $G$.  Let $\mathcal P = \bigcup_{u,v \in V} \mathcal P_{uv}$.
Then a {\em flow in $G$} is a mapping $F : \mathcal P \to [0,\infty)$.
For any $u,v \in V$, let $F[u,v] = \sum_{p \in \mathcal P_{uv}} F(p)$
be the amount of flow sent between $u$ and $v$.

Our main technical theorem concerns a class of flows we call {\em
  subset flows}.  Let $\mu$ be a probability distribution on subsets
of $V$.  Then $F$ is a {\em $\mu$-flow} if it satisfies
$F[u,v] = \Pr_{S \sim \mu}[u, v \in S]$ for all $u,v \in V$.
For $r \leq |V|$, we write $\mathcal F_r(G)$ for the set of all
$\mu$-flows in $G$ with $\supp(\mu) \subseteq {V \choose r}$.

We say a flow $F$ is an {\em integral flow} if it is supported on only
one path $p$ in each $\mathcal P_{uv}$, and a {\em unit flow} if $F[u,v] \in \{0,1\}$ for every
$u,v \in V$.
An edge-weighted graph $H$ is one which comes equipped
with a non-negative weight function $w : E(H) \to [0,\infty)$ on edges.
We say that a flow $F$ in $G$ is an {\em $H$-flow} if
there exists an injective mapping $\phi : V(H) \to V$ such that
for all $\{u,v\} \in E(H)$, we have $F[\phi(u),\phi(v)] \geq w(u,v)$.
In this case, $H$ is referred to as the {\em demand graph} and $G$ as
the {\em host graph.}

\medskip

We define the {\em squared $\ell_2$-congestion}, or simply {\em
  congestion}, of a flow $F$ by
$\vcon(F) = \sum_{v \in V} C_F(v)^2,$
where $C_F(v) = \sum_{p \in \mathcal P : v \in p} F(p)$.  This
congestion can also be written as
$$\vcon(F) = \sum_{p, p' \in \P} \sum_{v \in p \cap p'} F(p) F(p')$$
and is therefore bounded below by a more restricted
sum, the {\em intersection number}:
$$\inter(F)
   = \sum_{\substack{u,v,u',v' \\ |\{u,v,u',v'\}| = 4}}
        \sum_{\substack{p \in \P_{uv} \\ p' \in \P_{u'v'}}}
           \sum_{x \in p \cap p'} F(p) F(p').
$$

\section{Eigenvalues and spreading weights}
\label{sec:eigenv-spre-weights}

We will now reduce the problem of proving upper bounds on the eigenvalues of
a graph $G$, to the problem of proving lower lower bounds on the
congestion of subset flows in $G$.
In the present section, if $(X,d)$ is a metric space, and $x \in X, R \geq 0$,
we will use the notation
$$
B(x,R) = \{ y \in X : d(x,y) \leq R \}.
$$

\subsection{Padded partitions}
Let $(X,d)$ be a finite metric space.  We will view a partition
$P$ of $X$ as a collection of subsets, and also as a function
$P : X \to 2^X$ mapping a point to the subset that contains it.
We write $\beta(P,\Delta)$ for the infimal value of $\beta \geq 1$ such that
$$
\left|\left\{ \vphantom{\bigoplus} x \in X : B(x, \Delta/\beta) \subseteq P(x)\right\}\right|
\geq \frac{|X|}{2}.
$$
Let $\mathcal P_{\Delta}$ be the set of all partitions $P$
such that for every $S \in P$, $\diam(S) \leq \Delta$.
Finally, we define
$$\beta_{\Delta}(X,d) = \inf \left\{ \vphantom{\bigoplus}
\beta(P,\Delta) : P \in \mathcal P_{\Delta} \right\}.$$

The following theorem is a consequence \cite{Rao} of the main theorem
of Klein, Plotkin, and Rao \cite{KPR}, with the dependence of $r^2$
due to \cite{FT03}.

\begin{theorem}\label{thm:KPR}
  Let $G = (V, E)$ be a graph without a $K_{r,r}$ minor and $(V, d)$
  be any shortest-path semimetric on $G$, and let $\Delta > 0$.  Then
  $\beta_\Delta(V, d) = O(r^2)$.
\end{theorem}

In particular, if $G$ is planar then $\beta_\Delta(V, d)$ is upper bounded by
an absolute constant, and if $G$ is of genus $g > 0$ then
$\beta_\Delta(V, d) = O(g),$ since the
genus of $K_h$ is $\Omega(h^2)$ \cite[p.~118]{Harary94}.
The paper \cite{LSgenus} proves the following strengthening
(which is tight, up to a universal constant).

\begin{theorem}\label{thm:LSgenus}
Let $G=(V,E)$ be a graph of orientable genus $g$,
and $(V, d)$  be any shortest-path semimetric on $G$, and let $\Delta > 0$.  Then
  $\beta_\Delta(V, d) = O(\log g)$.
\end{theorem}

\subsection{Spreading vertex weights}

Consider
a non-negative weight function $\omega : V \to \mathbb R_+$ on vertices,
and extend $\omega$ to subsets $S \subseteq V$ via $\omega(S) = \sum_{v \in V} \omega(v)$.
We associate a vertex-weighted shortest-path metric by defining
$$
\dist_{\omega}(u,v) = \min_{p \in \mathcal P_{uv}} \omega(p).
$$
Say that $\omega$ is {\em $(r,\varepsilon)$}-spreading
if, for every $S \subseteq V$ with $|S| = r$, we have
$$
\frac{1}{|S|^2} \sum_{u,v \in S} \dist_{\omega}(u,v) \geq \varepsilon \sqrt{\sum_{v \in V} \omega(v)^2}.
$$

Write $\varepsilon_r(G, \omega)$ for the maximal value of $\varepsilon$
for which $\omega$ is $(r,\varepsilon)$-spreading.

\begin{theorem}[Higher eigenvalues]\label{thm:eigenvalues}
Let $G = (V,E)$ be any $n$-vertex
graph with maximum degree $d_{\max}$, and let $\lambda_k$ be the $k$th Laplacian
eigenvalue of $G$.  For any $k \geq 1$, the following holds.
For any weight function $\omega : V \to \mathbb R_+$ with
\begin{equation}\label{eq:omeganorm}
\sum_{v \in V} \omega(v)^2 = 1,
\end{equation}
we have
$$\lambda_{k} \lesssim \frac{d_{\max}}{\varepsilon^2n} \beta^2,$$
as long as
\begin{equation}\label{eq:newrestrict}
\frac{\beta^2}{\e^2} \leq \frac{n}{64},
\end{equation}
where $\varepsilon = \varepsilon_{\lfloor n/8k\rfloor}(G,\omega)$ and
$\beta = \beta_{\varepsilon/2}(V,\dist_{\omega})$.
\end{theorem}

\begin{proof}
Let $\omega$ be an $(\lfloor n/8k \rfloor, \varepsilon)$-spreading weight function.
Let $V = C_1 \cup C_2 \cup \cdots \cup C_m$ be a partition of $V$ into
sets of diameter at most $\varepsilon/2$, and define for every $i \in [m]$,
\begin{equation*}
\hat C_i  =\left\{\vphantom{\bigoplus} x \in C_i : B(x, \varepsilon/(2\beta)) \subseteq C_i \right\},
\end{equation*}
where $\beta = \beta_{\varepsilon/2}(V, \dist_{\omega})$.  By the definition of $\beta$,
there exists a choice of $\{C_i\}$ with 
\begin{equation}\label{eq:Csize}
|\hat C_1 \cup \hat C_2 \cup \cdots \cup \hat C_m| \geq n/2.
\end{equation}

Now, for any set $A \subseteq V$ with  $\diam(A) \leq \varepsilon/2$, we see that
\begin{equation}\label{eq:stupidref}
\frac{1}{|A|^2} \sum_{u,v \in A} \dist_{\omega}(u,v) \leq \frac{\varepsilon}{2} = \frac{\varepsilon}{2} \sqrt{\sum_{v \in V} \omega(v)^2}.
\end{equation}
Since $\diam(C_i) \leq \varepsilon/2$, if $|C_i| > n/8k$, then we could pass to a subset of $C_i$ of size
exactly $\lfloor n/8k\rfloor$ which satisfies \eqref{eq:stupidref}, but this would violate
the $(\lfloor n/8k \rfloor, \varepsilon)$-spreading property of $\omega$.  Hence we know that $|C_i| \leq n/8k$
for each $i=1,2,\ldots, m$.

Now, define a set of heavy vertices $H = \{ v \in V : \omega(v) \geq \frac{\e}{4\beta} \}$.  Clearly \eqref{eq:omeganorm} implies
that $|H| \leq 16\beta^2/\e^2 \leq n/4$, by assumption \eqref{eq:newrestrict}.
Finally, define for every $i=1,2,\ldots, m$, $\widetilde{C_i} = \hat C_i \setminus H$.
From \eqref{eq:Csize}, we know that $|\widetilde{C_1} \cup \widetilde{C_2} \cup \cdots \cup \widetilde{C_m}| \geq n/4$.

Thus by taking disjoint unions of the sets $\{\hat{C_i}\}$ which are each of size
at most $n/4k$, we can find sets $S_1, S_2, \ldots, S_{2k}$ with
\begin{equation}\label{eq:sizes}
\frac{n}{4k} \geq |S_i \setminus H| \geq \frac{n}{8k}.
\end{equation}
For each $i \in [2k]$, let $\tilde S_i$ be the $\e/(4\beta)$-neighborhood of $S_i$.
Observe that the sets $\{\tilde S_i\}$ are pairwise disjoint, since
by construction each is contained in a union of $C_i$'s, which are themselves pairwise disjoint.
Furthermore, by construction, we have $\dist_{\omega}(\tilde S_i, \tilde S_j) \geq \e/(2\beta)$ for $i \neq j$.

\remove{
Additionally, we have the estimate, for every $i \in [2k]$,
\begin{equation}\label{eq:ubound}
|\tilde S_i| \leq n - \sum_{j \neq i} |S_j| \leq n - (2k-1) \frac{n}{4k} \leq \frac{3n}{4}.
\end{equation}
}
Now define, for every $i \in [2k]$, define
$$
W(\tilde S_i) = \sum_{u \in \tilde S_i} \sum_{v : uv \in E} \left[\omega(u)+\omega(v)\right]^2
$$
Clearly, we have
$$
\sum_{i=1}^{2k} W(\tilde S_i) \leq 2 \sum_{uv \in E} [\omega(u)+\omega(v)]^2 \leq 4 d_{\max} \sum_{v \in V} \omega(v)^2 = 4 d_{\max},
$$
where the latter equality is \eqref{eq:omeganorm}.
Hence if we renumber the sets so that
$\left\{\tilde S_1, \tilde S_2, \ldots, \tilde S_k\right\}$ have the smallest $W(\tilde S_i)$ values, then
for each $i=1,2,\ldots, k$, we have
$W(\tilde S_i) \leq \frac{4 d_{\max}}{k}$.

Finally, we define functions $f_1, f_2, \ldots, f_k : V \to \mathbb R$
by
$$f_i(x) = \begin{cases}
0 & \textrm{if } \omega(x) \geq \frac{\e}{4\beta}, \\
\max\left\{0, \frac{\varepsilon}{4\beta} - \dist_{\omega}(x,S_i)\right\} & \textrm{otherwise.}
\end{cases}
$$
so that $f_i$ is supported on $\tilde S_i$,
and $f_i(x) = \varepsilon/(4\beta)$ for $x \in S_i \setminus H$.

Observe that $f_i|_{V \setminus H}$ is $1$-Lipschitz, hence
\begin{eqnarray*}
\|f_i\|_{\mathcal L}^2 = \sum_{uv \in E} |f_i(u)-f_i(v)|^2 &\leq&
\sum_{u \in \tilde S_i \cap H} \deg(u) \left(\frac{\e}{4\beta}\right)^2
+
 \sum_{u \in \tilde S_i \setminus H} \mathop{\sum_{v : uv \in E}}_{v \notin H} |f_i(u)-f_i(v)|^2 \\
&\leq & \sum_{u \in \tilde S_i \cap H} \deg(u)\, \omega(u)^2 + \sum_{u \in \tilde S_i \setminus H} \mathop{\sum_{v : uv \in E}}_{v \notin H} \dist_{\omega}(u,v)^2  \\
&\leq & \sum_{u \in \tilde S_i} \sum_{v : uv \in E} \left[\omega(u)+\omega(v)\right]^2 \\
&=& W(\tilde S_i) \leq \frac{4 d_{\max}}{k}.
\end{eqnarray*}
Furthermore the functions $\{f_i\}$ satisfy,
$$
\|f_i\|^2\geq \left(\frac{\varepsilon}{4\beta}\right)^2 |S_i \setminus H|
  \geq \frac{\varepsilon^2}{128 \beta^2} \frac{n}{k},
$$
where in the final inequality we have used \eqref{eq:sizes}.

Combining the preceding two estimates shows that for each $f_i$,
$$
\frac{\|f_i\|_{\mathcal L}^2}{\|f_i\|^2}
  \leq \frac{d_{\max}}{128 n} \left(\frac{\beta}{\varepsilon}\right)^2.
$$

Thus by Lemma \ref{lem:simple}, the proof is complete once we show that $\supp(f_i) \cap N(\supp(f_j)) = \emptyset$ when $i \neq j$.
To this end, consider a pair $x,y \in V$ with $x \sim y$, but $x \in \tilde S_i$ and $y \in \tilde S_j$ for $i \neq j$.
Since $\dist_w(x,y) \geq \dist(\tilde S_i, \tilde S_j) \geq \e/(2\beta)$, it must be that either $\omega(x) \geq \e/(4\beta)$ or $\omega(y) \geq \e/(4\beta)$,
implying that either $f_i(x)=f_j(x)=0$ or $f_i(y)=f_j(y)=0$.
\end{proof}

\remove{
The following results follow from Theorem~\ref{thm:eigenvalues},
Theorem~\ref{thm:KPR} and Theorem~\ref{thm:LSgenus}, Theorem~\ref{thm:duality}, Lemma~\ref{lem:round},
Theorem~\ref{thm:subset-crossings}, and Corollary~\ref{cor:ixn-basic}.
\begin{corollary}
\label{cor:specificbounds}
  If $G$ is planar, then
  $$ \lambda_k \leq O\left( d_{\max} \frac k n \right). $$
  If $G$ is of genus $g > 0$, then
  $$ \lambda_k \leq O\left( d_{\max}\frac k n   g (\log g)^2 \right). $$
  If $G$ is $K_h$-minor-free, then
  $$ \lambda_k \leq O\left( d_{\max}\frac k n h^6 \log h \right). $$
\end{corollary}
}

\subsection{Spreading weights and subset flows}
\label{sec:spre-weights-subset-flows}

We now show a duality between the optimization problem of
finding a spreading weight $\omega$ and the problem of minimizing
congestion in subset flows.
The following theorem is proved by a standard
Lagrange multipliers argument.

\begin{theorem}[Duality]\label{thm:duality}
  Let $G = (V, E)$ be a graph and let $r \leq |V|$.
  Then
  $$
  \max \left\{ \vphantom{\bigoplus} \varepsilon_r(G,\omega)
              \middle| \omega : V \to \mathbb R_+ \right\}
  = \frac{1}{r^2} \min \left\{ \vphantom{\bigoplus} \sqrt{\vcon(F)}
              \middle| F \in \mathcal F_r(G) \right\}
  .
  $$
\end{theorem}

\begin{proof}
  We shall write out the optimizations $\max_\omega \varepsilon_r(G,
  \omega)$ and $\frac{1}{r^2}\min_F \sqrt{\vcon(F)}$ as convex programs, and show that
  they are dual to each other.  The equality then follows from
  Slater's condition \cite[Ch. 5]{convex-optimization}:
  \begin{fact}[Slater's condition for strong duality]
    When the feasible region for a convex program $(\mathbf P)$ has
    non-empty interior, the values of $(\mathbf P)$ and its dual
    $(\mathbf P^*)$ are equal.
  \end{fact}

  We begin by expanding $\max_\omega \varepsilon_r(G, \omega)$ as a
  convex program $(\mathbf P)$.  Let $P \in \{0, 1\}^{\mathcal P
    \times V}$ be the path incidence matrix; $Q \in \{0, 1\}^{\mathcal
    P \times \binom V 2}$ the path connection matrix; and $R \in \{0,
  1\}^{\binom V r \times \binom V 2}$ a normalized set containment
  matrix, respectively defined as
  $$ P_{p, v} = \pieceww 1 {v \in p} 0 {\textrm{else}}
  \qquad \qquad
  Q_{p, uv} = \pieceww 1 {p \in \mathcal P_{uv}} 0 {\textrm{else}}
  \qquad \qquad
  R_{S, uv} = \pieceww {1/r^2} {\{u, v\} \subset S} 0 {\textrm{else.}}
  $$
  Then the convex program $(\mathbf P) = \max_\omega \varepsilon_r(G,
  \omega)$ is
  $$
  \begin{array}{rccc}
    \textrm{minimize} & -\varepsilon \\
    \textrm{subject to} &
      \varepsilon \mathbf 1 \preceq R d &
      Q d \preceq P s &
      s^\top s \leq 1 \\
    & d \succeq 0 & s \succeq 0
  \end{array} .
  $$
  Introducing the non-negative Lagrange multipliers $\mu, \lambda, \nu$,
  the Lagrangian function is
  $$ L(d, s, \mu, \lambda, \nu) =
   -\varepsilon +
   \mu^\top (\varepsilon \mathbf 1 - R d) +
   \lambda^\top (Q d - P s) +
   \nu (s^\top s - 1)
  $$
  so that $(\mathbf P)$ and its dual $(\mathbf P^*)$ may be written as
  \begin{align*}
  (\mathbf P)   &= \inf_{\varepsilon, d, s} \sup_{\mu, \lambda, \nu}
                     L(d, s, \mu, \lambda, \nu) \\
  (\mathbf P^*) &= \sup_{\mu, \lambda, \nu} \inf_{\varepsilon, d, s}
                     L(d, s, \mu, \lambda, \nu).
  \end{align*}
  Now we simplify $(\mathbf P^*)$.  Rearranging terms in $L$, we have
  \begin{align*}
    (\mathbf P^*) &=
      \sup_{\mu, \lambda, \nu} \inf_{\varepsilon, d, s}
        (\mu^\top \mathbf 1 - 1) \varepsilon +
        (\lambda^\top Q - \mu^\top R) d +
        (\nu s^\top s - \lambda^\top P s)
        - \nu \\
    &= \sup_{\mu, \lambda, \nu}
        \inf_\varepsilon (\mu^\top \mathbf 1 - 1) \varepsilon +
        \inf_d (\lambda^\top Q - \mu^\top R) d +
        \inf_s (\nu s^\top s - \lambda^\top P s)
        - \nu.
  \end{align*}
  Now the infima $\inf_\varepsilon (\mu^\top \mathbf 1 - 1)
  \varepsilon$ and $\inf_d (\lambda^\top Q - \mu^\top R) d$ are either
  0 or $-\infty$, so at the optimum they must be zero and $\mu^\top
  \mathbf 1 - 1 \geq 0$, $\lambda^\top Q - \mu^\top R \succeq 0$.
  With these two constraints, the optimization reduces to
  $\sup_{\lambda, \nu} \inf_s (\nu s^\top s - \lambda^\top P s) - \nu$.  At the
  optimum, the gradient of the infimand is zero, so $s = \frac{P^\top
    \lambda}{2\nu}$ and the infimum is $-\frac{\|P^\top \lambda\|_2^2}{4\nu}$.
  Then at the maximum,
  $\nu = \frac12 \|P^\top \lambda\|_2$, so that the
  supremand is $- \|P^\top \lambda\|_2$.  We have shown that $(\mathbf
  P^*)$ is the convex program
  $$
  \begin{array}{rcc}
    \textrm{maximize} & -\left\|P^\top \lambda\right\|_2 \\
    \textrm{subject to} &
      \lambda^\top Q \succeq \mu^\top R &
      \mu^\top \mathbf 1 \geq 1 \\
    & \lambda \succeq 0 & \mu \succeq 0
  \end{array} .
  $$
  This program is precisely (the negative of) the program to minimize
  vertex 2-congestion of a subset flow in $\mathcal F_r(G)$, where the
  subset weights are normalized to have unit sum.  The proof is
  complete.
\end{proof}

\remove{
\begin{theorem}\label{thm:spreading-weights}
  Let $G = (V, E)$ be a graph and let $r \leq |V|$.  Then there is a
  vertex weighting $\omega : V \to \mathbb R_+$ such that
  $$ \varepsilon_r(G,\omega)
     \geq \sqrt{\frac{r}{c(\inter_G) n}}. $$
\end{theorem}
\begin{proof}
  Observe that $\vcon_2^G G \geq \sqrt{\inter_G(G)}$, and combine
  Theorem~\ref{thm:duality} and Corollary~\ref{cor:subset-crossings}.
\end{proof}
}

\section{Congestion measures}
\label{sec:cong-meas}

In this section, we develop concepts that will enable us to give lower
bounds on the congestion $\vcon(F)$ of all subset flows $F$ in a given
graph $G$.  The reader may wish to consult with Section \ref{sec:prelim-flows}
to recall the relevant definitions.

\begin{definition}
  Let $G$ be an arbitrary host graph, and $H$ an edge-weighted demand graph.
Define the \emph{$G$-congestion
    of $H$} by
  $$ \vcon_G(H)
   = \min_{\textrm{\rm $F$ an $H$-flow in $G$}}  \vcon(F)
  $$
  and the \emph{$G$-intersection number of $H$} by
  $$ \inter_G(H) = \min_{\textrm{\rm $F$ an $H$-flow in $G$}} \inter(F),$$
  and the \emph{integral $G$-intersection number of $H$} by
  $$ \inter^*_G(H) = \min_{\textrm{\rm $F$ an integral $H$-flow in $G$}} \inter(F).$$
\end{definition}

Note that even if $H$ is a unit-weighted graph and $\inter^*_G(H)=0$,
this does not imply that $G$ contains an $H$-minor.  This is because
the intersection number involves quadruples of four distinct vertices.
For example, if $H$ is a triangle, then $\inter^*_G(H)=0$ for any $G$,
even when $G$ is a tree (and thus does not have a triangle as a minor).
However, we recall the following (which appears as Lemma 3.2 in \cite{BLR}).

\begin{lemma}[\cite{BLR}]\label{lem:interminor}
If $H$ is a unit-weighted, bipartite demand graph in which every node
has degree at least two, then for any graph $G$, $\inter^*_G(H)=0$ implies
that $G$ contains an $H$-minor.
\end{lemma}

The next lemma is proved via randomized rounding.

\begin{lemma}[Rounding]\label{lem:round}
  For any graph $G$ and unit flow $F$, there is an integral unit flow
  $F^*$ with $F^*[u,v] = F[u,v]$ for all $u,v \in V(G)$, and such that
  $$ \inter(F^*) \leq \inter(F) . $$
  Consequently for every $G$ and unit-weighted $H$,
  \begin{equation} \label{eq:round}
    \inter^*_G(H) = \inter_G(H) \leq \vcon_G(H).
  \end{equation}
\end{lemma}

\begin{proof}
  We produce an integral flow $F^*$ randomly by rounding $F$.  For
  each pair of endpoints $u,v$, choose independently a path $p_{uv}$ in
  $\P_{uv}$ with $\Pr[p_{uv} = p] = F(p)$ for each $p$.  Then
  \begin{multline*}
  \E[\inter(F^*)]
   = \sum_{\substack{u,v,u',v' \\ |\{u,v,u',v'\}| = 4}}
           \E\big[\left|p_{uv} \cap p_{u'v'}\right|\big]
   = \sum_{\substack{u,v,u',v' \\ |\{u,v,u',v'\}| = 4}} \,
        \sum_{\substack{p \in \P_{uv} \\ p' \in \P_{u'v'}}}
           \sum_{x \in p \cap p'} F(p) F(p')
   = \inter(F)
  \end{multline*}
  so that with positive probability we must have $\inter(F^*) \leq
  \inter(F)$.  Equation~\eqref{eq:round} follows because $\inter(F) \leq
  \vcon(F)$ always.
\end{proof}

\begin{definition}\label{def:conmeasure}
Given a host graph $G$, we say that $\inter_G$ is a {\em $(c,a)$-congestion measure}
if for all unit-weighted graphs $H=(V,E)$, we have the inequality
\begin{equation}
\label{eq:E3V2}
\inter_G(H) \geq \frac{|E|^3}{c |V^2|} - a |V|\,.
\end{equation}
  In particular, $\mathcal \inter_G(K_n) \geq \frac{n^4}{8c} - an$.
\end{definition}

\remove{The term $a(\mathcal C)|V|$ in (\ref{eq:E3V2}) is present for
technical reasons.
Throughout this and the following section, terms containing
$a(\mathcal C)$ should be thought of as negligible and ignored on a
first reading.}

\begin{lemma}\label{lem:inter}
  Suppose that for some $G$ and $k = k(G)$, every unit-weighted $H$ obeys
  \begin{equation}
    \label{eq:inter-weak}
    \inter^*_G(H) = \inter_G(H) \geq |E(H)| - k |V(H)| - k^2.
  \end{equation}
  Then it follows that for every unit-weighted $H$,
  \begin{equation}
    \label{eq:inter-strong}
    \inter_G(H) \geq \frac 1 {27} \frac{|E(H)|^3}{k^2 |V(H)|^2}
                   - k |V(H)|
  \end{equation}
  so that $\inter_G$ is an $(27 k^2, k)$-congestion measure.
\end{lemma}
\begin{proof}
  It suffices to consider $|E(H)| \geq 3 k |V(H)|$ since
  otherwise the right-hand side of inequality (\ref{eq:inter-strong})
  is negative.

  Fix any $H$-flow $F$ in $G$.  Sample the nodes of $H$
  independently with probability $p$ each to produce a new demand
  graph $H'$ and flow $F' = F|_{H'}.$  Then $\inter(F') \geq \inter_G(H') \geq
  |E(H')| - k |V(H')| - k^2$, and by taking expectations we have
  $$ p^4 \inter(F) \geq p^2 |E(H)| - p k |V(H)| - k^2.$$
  Choosing $p = 3 k |V(H)|/|E(H)|$ and using the fact that
  $|E(H)|/|V(H)|^2 < 1$ we obtain (\ref{eq:inter-strong}).
\end{proof}

The next proof follows employs the techniques of \cite{BLR}.

\begin{corollary}\label{cor:ixn-basic}
  If $G$ is planar, then $\inter_G$ is an $(O(1),3)$-congestion measure.
  If $G$ has genus $g > 0$, then $\inter_G$ is an $(O(g), O(\sqrt{g}))$-congestion measure.
  If $G$ is $K_h$-minor-free, then $\inter_G$ is an $(O(h^2 \log h), O(h \sqrt{\log h}))$-congestion measure.
\end{corollary}

\begin{proof}
Suppose that $H$ is a unit-weighted demand graph.
  If $F$ is an integral $H$-flow with $\inter(F) > 0$, then some path
  in $F$ and corresponding edge of $H$ can be removed to yield an
  integral $H'$-flow $F'$ with $\inter(F') \leq \inter(F) - 1$.
  Therefore to prove (\ref{eq:inter-weak}) it suffices to consider $H$
  with $\inter_G(H) = 0$ and show that $|E(H)| \leq k |V(H)| + k^2$.
  Then Lemma~\ref{lem:inter} will imply $\inter_G$ is an
  $(O(k^2), k)$-congestion measure.

  When $G$ is planar, an $H$-flow $F$ in $G$ with $\inter(F) = 0$
  gives a drawing of $H$ in the plane without crossings, so that $H$
  itself is planar.  Then an elementary application of the Euler
  characteristic gives
  $$ |E(H)| \leq 3|V(H)| - 6 < 3|V(H)|. $$

  When $G$ is of genus at most $g > 0$, the same argument gives
  $$ |E(H)| \leq 3|V(H)| + 6(g-1),$$
  which suffices for $k = O(\sqrt g)$.

  For $K_h$-minor-free $G$ and $H$ with $\inter_G(H) = 0$, if $H$ is
  bipartite with minimum degree 2, then Lemma \ref{lem:interminor} implies that $H$ is
  $K_h$-minor-free, so that $|E(H)| \leq c_{KT} |V(H)| h \sqrt{\log h}$
  by the theorem of Kostochka \cite{Kostochka} and Thomason
  \cite{Thomason}.

  For general $H$, we can first take a partition to obtain
  a bipartite subgraph $H'$ with $|E(H')| \geq |E(H)|/2$.
  We then remove isolated vertices from $H'$, and iteratively
  remove vertices of degree one and the associated edges
  to obtain a bipartite subgraph $H''$ with minimum degree two, and
  \begin{equation}\label{eq:intermed1}
  |E(H'')| \geq |E(H')| - |V(H')| \geq |E(H)|/2 - |V(H)|\,.
  \end{equation}
  Now $\inter_G(H)=0$ together with Lemma \ref{lem:interminor} implies that
  $$
  |E(H'')| \leq 2 c_{KT}h \sqrt{\log h} |V(H'')|\,
  $$
  which together with \eqref{eq:intermed1}, implies that
  $|E(H)| \leq O(h \sqrt{\log h}) |V(H)|$.
  \end{proof}

In the next section, we will also require the following lemma.

\begin{lemma}\label{lem:SS'4}
  Let $\mu$ be any probability distribution over subsets of $V$.  Writing $H_\mu$
  for the graph on $V$ with edge weights $H_\mu(u, v) = \Pr_{S \sim
    \mu} [u, v \in S]$, we have
  $$ \inter_G(H_\mu) \geq \mathbb E_{S \sim \mu, S' \sim \mu} \left[\inter_G(K_{|S \cap S'|})\right],$$
  where by $K_n$ we intend the unit-weighted complete graph on $n$ vertices.
\end{lemma}

\begin{proof}
  Let $F$ be any $H_{\mu}$-flow, and let the vertices of $H_\mu$ be identified with the
  corresponding vertices of $G$ (recall that every $H$-flow in $G$ comes with
  an injection from $V(H)$ into $V(G)$).

  Now, for every $u,v \in V$, $S \subseteq V$, and $p \in \mathcal P_{uv}$, define the flow $F^S$ by,
  $$
  F^S(p) = \begin{cases} \displaystyle \frac{F(p)}{F[u,v]} \mu(\{S\}) & u,v \in S \textrm{ and } F[u,v] \neq 0\\
  0 & \textrm{otherwise,}
  \end{cases}
  $$
  and observe that
  since $F[u,v] = \Pr_{S \sim \mu} [u,v \in S]$, we have $F = \sum_{S \subseteq V} F^S$.

  In this case, we can write
  \begin{eqnarray*}
  \inter(F) &=& \sum_{\substack{u,v,u',v' \\ |\{u,v,u',v'\}| = 4}}
          \sum_{\substack{p \in \P_{uv} \\ p' \in \P_{u'v'}}} |p \cap p'|\,F(p) F(p') \\
  &=&
  \sum_{\substack{u,v,u',v' \\ |\{u,v,u',v'\}| = 4}}
          \sum_{\substack{p \in \P_{uv} \\ p' \in \P_{u'v'}}} |p \cap p'|\,\left(\sum_{S} F^S(p)\right) \left(\sum_{S}  F^S(p')\right) \\
  &=&
    \sum_{\substack{u,v,u',v' \\ |\{u,v,u',v'\}| = 4}}
          \sum_{\substack{p \in \P_{uv} \\ p' \in \P_{u'v'}}} |p \cap p'| \sum_{S,S'} F^S(p) F^{S'}(p') \\
  &= &
  \sum_{S,S'} \sum_{\substack{u,v \in S \\ u',v' \in S' \\ |\{u,v,u',v'\}|=4}} \sum_{\substack{p \in \P_{uv} \\ p' \in \P_{u'v'}}} |p \cap p'| F^S(p) F^{S'}(p') \\
  &=&
  \sum_{S,S'} \sum_{\substack{u,v \in S \\ u',v' \in S' \\ |\{u,v,u',v'\}|=4}} \sum_{\substack{p \in \P_{uv} \\ p' \in \P_{u'v'}}} |p \cap p'|
   \mu(\{S\}) \frac{F(p)}{F[u,v]} \mu(\{S'\}) \frac{F(p')}{F[u',v']}\\
   &\geq&
  \sum_{S,S'} \mu(\{S\}) \mu(\{S'\}) \sum_{\substack{u,v,u',v' \in S \cap S' \\ |\{u,v,u',v'\}|=4}} \sum_{\substack{p \in \P_{uv} \\ p' \in \P_{u'v'}}} |p \cap p'|
   \frac{F(p)}{F[u,v]} \frac{F(p')}{F[u',v']}\\
&\geq &
\mathbb \E_{S \sim \mu, S' \sim \mu} \left[\inter_G(K_{|S \cap S'|})\right]\,,
     \end{eqnarray*}
where we have used the fact that the double sum in the penultimate line contains precisely
the intersection number of a unit-weighted complete-graph flow on $S \cap S'$.
\end{proof}

\section{Congestion for subset flows}
\label{sec:subset-congestion}

We now prove our main estimate on
the congestion incurred by subset flows
in terms of a graph's congestion measure.

\medskip

The proof of Theorem \ref{thm:subset-crossings} below involves
some delicate combinatorial and analytic arguments.
The difficulty lies in controlling
the extent to which $\mu$ is a mixture of three different types of ``extremal''
distributions:
\begin{enumerate}
\item $\mu$ is uniformly distributed on all sets of size $r$,
\item $\mu$ is concentrated on a single set of size $r$,
\item $\mu$ is uniform over $n/r$ disjoint sets, each of size $r$.
\end{enumerate}

In the actual proof, we deal with the corresponding cases:
(1') $\mu$ is ``uniformly spread'' over edges, i.e.
$\Pr_{S \sim \mu} [u,v \in S]$ is somewhat uniform over
choices of $u,v \in V$.
In this case, we have to take a global approach, showing that not
only are there many intra-set crossings, but also a lot of crossing
weight is induced by crossing edges coming from different sets.
(2') $\Pr_{S \sim \mu} [u \in S]$ is unusually large for all $u \in V'$
with $|V'| \ll |V|$.  In this case, there is a ``density increment''
on the induced subgraph $G[V']$, and we can apply induction.
Finally, if we are in neither of the cases (1') or (2'), we are left
to show that, in some sense, the distribution $\mu$ must
be similar to case (3) above, in which case we can appeal
to the classical dense crossing bounds applied to the complete graph on
$S \cap S'$ where $S,S' \sim \mu$ are chosen i.i.d.

\begin{theorem}\label{thm:subset-crossings}
There is a universal constant $c_0 > 0$ such that the following holds.
Let $\mu$ be any probability distribution on subsets of $[n]$.
For $u,v \in [n]$, define $$F(u,v) = \Pr_{S \sim \mu}[u,v \in S]$$
and let $H_\mu$ be the graph on $[n]$ weighted by $F$.
For any graph $G$ such that $\inter_G$ is a $(c,a)$-congestion measure, we have
$$
\inter_G(H_\mu) \gtrsim \frac{1}{c n}\left(\E|S|^2\right)^{5/2}
            \; - c_0 \frac{a}n \E|S|^2
$$
\end{theorem}

\begin{corollary}\label{cor:subset-crossings}
  If $\mu$ is supported on ${[n] \choose r}$ for some $r$, then
  $ \inter_G(H_\mu) \gtrsim \frac{r^5}{c n} - c_0 \frac{a r^2}n. $
  In particular, if $r \gtrsim (a\cdot c)^{1/3}$, then
  $$ \inter_G(H_\mu) \gtrsim \frac{r^5}{c n}. $$
\end{corollary}

\medskip

\begin{proof}[Proof of Theorem~\ref{thm:subset-crossings}]
We will freely use the fact that
$$
\mathbb E|S|^2 = \sum_{u,v} F(u,v).
$$
Also, put
$F(u) = \Pr_{S \sim \mu} [u \in S]$ for $u \in [n]$.

\medskip

The proof will proceed by induction on $n$, and
will be broken into three cases.
Let
$$\beta = \sqrt{\frac{1}{n^2} \sum_{u,v} F(u,v)},$$
and put $E(\alpha',\alpha) =  \{(u,v) : \alpha' \leq F(u,v) \leq \alpha \}$.
Define the set of ``heavy vertices'' as
$$
H_{K} = \{ u : F(u) \geq K \beta \},
$$
for some constant $K \geq 1$ to be chosen later.
Let $E_H = \{ (u,v) : u,v \in H_K \}$ and $E_{HL} = \overline{E(0,\beta) \cup E_H}$.

\medskip

\noindent
{\bf Case I (Light edges):} $\displaystyle\sum_{(u,v) \in E(0,\beta)} F(u,v) \geq \frac14 \sum_{u,v} F(u,v).$

\medskip

The desired conclusion comes from applying the following claim.


\begin{claim}
\label{claim:lite}
For every $\beta \in [0,1]$, we have
\begin{equation}\label{eq:lite}
\inter_G(H_\mu) \gtrsim
   \frac{\left(\sum_{(u,v) \in E(0,\beta)} F(u,v)\right)^3}
        {\beta c n^2}
   - 2 \beta^2 a n.
\end{equation}
\end{claim}

\begin{proof}
First, observe that by (\ref{eq:E3V2}), the subgraph consisting of the
edges in $E(\alpha, \beta)$ contributes at least
$$ \alpha^2 \frac{|E(\alpha,\beta)|^3}{c n^2}
   - \beta^2 a n $$
to $\inter_G(H_\mu)$ for every $\alpha, \beta \in [0,1]$.  Therefore
letting $E_i = E\left(2^{-i-1} \beta, 2^{-i} \beta\right)$, we have
$$
\inter_G(H_\mu) \gtrsim \frac{1}{c n^2}
                  \sum_{i=0}^{\infty} 2^{-2i} \beta^2 |E_i|^3
                - a n \sum_{i=0}^{\infty} 2^{-2i} \beta^2 .
$$
Let $F_i = \sum_{(u,v) \in E_i} F(u,v)$ so that
$|E_i| \geq (2^i/\beta) F_i$, and then
$$
\inter_G(H_\mu) \gtrsim \frac{1}{\beta c n^2}
                  \sum_{i=0}^{\infty} 2^{i} F_i^3
                - 2 \beta^2 a n,
$$
but also $\sum_{i=0}^\infty F_i = \sum_{u,v \in E(0,\beta)} F(u,v)$.
Thus \eqref{eq:lite} is proved by noting that
$$
\sum_{i=0}^{\infty} F_i
= \sum_{i=0}^{\infty} \left(2^{-i/3} \cdot 2^{i/3} F_i\right)
\leq \left(\sum_{i=0}^{\infty} 2^{-i/2}\right)^{2/3}
     \left(\sum_{i=0}^\infty 2^i F_i^3\right)^{1/3}
< 2.27 \left(\sum_{i=0}^\infty 2^i F_i^3\right)^{1/3},
$$
using H\"older's inequality.
\end{proof}

\medskip

\noindent
{\bf Case II (Heavy endpoints): $\displaystyle\sum_{(u,v) \in E_H} F(u,v) \geq \frac14 \sum_{u,v} F(u,v)$.}

\medskip
Observe that
$$
\sum_{u \in [n]} F(u) = \mathbb E_{S \sim \mu} |S| \leq \sqrt{\mathbb E_{S \sim \mu} |S|^2} = \sqrt{\sum_{u,v} F(u,v)} = \beta n,
$$
hence $|H_{K}| \leq n/K$ by Markov's inequality.

Apply the statement
of the Theorem inductively to the distribution over
subsets of $V(H_K)$ corresponding to the
random set $S \cap V(H_K)$, to conclude that
\begin{equation}\label{eq:heavy1}
\inter_G(H_\mu) \gtrsim \frac{K}{c n}
                  \left( \sum_{(u,v) \in E_H} F(u,v) \right)^{5/2}
            \;  - c_0 \frac{a}n \sum_{(u,v) \in E_H} F(u,v)
.
\end{equation}
Consequently, by choosing $K=32$,
under the assumption of this case,
$$
K \left( \frac{\sum_{(u,v) \in E_H} F(u,v)}{\sum_{u,v} F(u,v)}\right)^{5/2} \geq 1
$$
and the conclusion again follows.

\bigskip
\noindent
{\bf Case III (Heavy edges, light endpoints):} $\displaystyle\sum_{(u,v) \in E_{HL}} F(u,v) \geq \frac12 \sum_{u,v} F(u,v)$.

\medskip
By definition,
$E_{HL} = \{ (u,v) : F(u,v) > \beta, \{u,v\} \nsubseteq H_K \}$.
Let $\kappa = (16 ac)^{1/3}$, so that $\frac{\kappa^4}{8c}  \geq 2 a \kappa$.
Using Lemma~\ref{lem:SS'4} and the fact that $\inter_G$ is a $(c,a)$-congestion measure, we have
\begin{align}
\inter_G(H_\mu) &\geq \E_{S \sim \mu, S' \sim \mu} \left[\inter_G(K_{|S \cap S'|})\right]\nonumber \\
&\geq
\E_{S \sim \mu, S' \sim \mu} \left[\inter_G(K_{|S \cap S'|}) \1_{|S \cap S'| \geq \kappa}\right] \nonumber \\
 &\geq \frac1{8c}
                  \mathbb E_{S \sim \mu, S' \sim \mu}\left[|S \cap S'|^4 \,\1_{|S \cap S'| \geq \kappa}\right]
        - a \mathbb E_{S \sim \mu, S' \sim \mu}\left[|S \cap S'| \,\1_{|S \cap S'| \geq \kappa}\right] \nonumber\\
        &\geq \frac1{16 c}
                  \mathbb E_{S \sim \mu, S' \sim \mu}\left[|S \cap S'|^4 \,\1_{|S \cap S'| \geq \kappa}\right] \nonumber \\
                  &=
\frac{1}{16c} \sum_{u \in [n]} \Pr[u \in S \cap S'] \,\mathbb E_{S \sim \mu, S' \sim \mu}\left[|S \cap S'|^3 \,\1_{|S \cap S'| \geq \kappa}\,\Big|\, u \in S \cap S'\right]  \nonumber \\
&=
\frac{1}{16c} \sum_{u \in [n]} (\Pr[u \in S])^2 \,\mathbb E_{S \sim \mu, S' \sim \mu}\left[|S \cap S'|^3 \,\1_{|S \cap S'| \geq \kappa}\,\Big|\, u \in S \cap S'\right]  \nonumber \\
&\geq
\frac{1}{16c} \sum_{u :\beta \leq F(u) \leq K \beta} F(u)^2 \,\mathbb E_{S \sim \mu, S' \sim \mu}\left[|S \cap S'|^3 \,\1_{|S \cap S'| \geq \kappa}\,\Big|\, u \in S \cap S'\right] \nonumber \\
&\geq
\frac{\beta^2}{16c} \sum_{u : \beta \leq F(u) \leq K \beta} \mathbb E_{S \sim \mu, S' \sim \mu}\left[|S \cap S'|^3\,\Big|\, u \in S \cap S'\right]
- \frac{K^2 \beta^2}{16c} n \kappa^3.
\nonumber
\remove{
\label{eq:heavylight1}
 &= \frac1{c}
    \sum_{u \in [n]}  \left(\Pr[u \in S]\right)^2
      \sum_{v,v',v'' \in [n]} \left(
        \Pr[v,v',v'' \in S \mid u \in S] \right)^2 \\
 &\phantom=
   - a \sum_{u \in [n]}  \left(\Pr[u \in S]\right)^2 \nonumber
.}
\end{align}
Since $\frac{K^2 \beta^2}{16 c} n \kappa^3 = \frac{K^2 a}{n} \mathbb E|S|^2$, to finish
the proof we need only show that
\begin{equation}\label{eq:final}
\sum_{u : \beta \leq F(u) \leq K \beta} \mathbb E_{S \sim \mu, S' \sim \mu}\left[|S \cap S'|^3\,\Big|\, u \in S \cap S'\right] \gtrsim n \left(\mathbb E |S|^2\right)^{3/2}.
\end{equation}

Now for each $u \in [n]$ with $F(u) = \Pr[u \in S] > 0$, let $\mu_u$ denote the distribution $\mu$
conditioned on $u \in S$.  Let $I_{vS}$ denote the
indicator of the event $\{v \in S\}$, so that $\Pr[v \in S \mid u \in S] = \E_{S \sim \mu_u}[I_{vS}]$.
In this case,
\begin{align}
\mathbb E_{S \sim \mu, S' \sim \mu}\left[|S \cap S'|^3\,\Big|\, u \in S \cap S'\right]
  &= \sum_{v,v',v'' \in [n]}
       \E_{S \sim \mu_u, S' \sim \mu_u}\big[
         I_{vS}I_{v'S}I_{v''S} I_{vS'}I_{v'S'}I_{v''S'}\big] \nonumber \\
  &= \E_{S \sim \mu_u, S' \sim \mu_u} \Big[
       \Big( \sum_v I_{vS}I_{vS'} \Big)^3
       \Big] \nonumber \\
&\geq \Big( \E_{S \sim \mu_u, S' \sim \mu_u} \Big[
       \sum_v I_{vS}I_{vS'}
       \Big] \Big)^3 \nonumber \\
  &= \Big( \sum_v \left(\E_{S \sim \mu_u} \big[I_{vS}\big]\right)^2
       \Big)^3\,. \nonumber
\end{align}
Therefore the left hand side of \eqref{eq:final} is at least
\begin{align}
          \sum_{u : \beta \leq F(u) \leq K \beta}
       \left( \sum_v \Pr[v \in S \mid u \in S]^2 \right)^3
 &\geq \frac1{K^6}
     \sum_{u : \beta \leq F(u) \leq K \beta}
       \left|\{ v : F(u,v)/F(u) \geq 1/K \}\right|^3
     \nonumber \\
 &\geq \frac1{K^6}
     \sum_{u : \beta \leq F(u) \leq K \beta}
       \left|\{ v : F(u, v) \geq \beta \}\right|^3
      \label{eq:heavylight3} \\
      &\geq \frac{1}{K^6}
     \sum_{u : \beta \leq F(u) \leq K \beta}
       \left|\{ v : (u,v) \in E_{HL} \}\right|^3, \nonumber
\end{align}
since each of the edges in $E_{HL}$ appears at least once in the sum \eqref{eq:heavylight3},
because every edge $(u,v) \in E_{HL}$ has either $F(u) \leq K\beta$ or
$F(v) \leq K\beta$.

In particular, for such edges, $F(u,v) \leq K\beta$, which means
that
\begin{equation}\label{eq:EHL}
|E_{HL}| \geq
 \frac{\sum_{(u,v) \in E_{HL}} F(u,v)}{K\beta}.
\end{equation}

Thus by the power-mean inequality,
the left hand side of \eqref{eq:final} is at least
\begin{align}
     \frac{1}{K^6}
     \sum_{u : \beta \leq F(u) \leq K \beta}
       \left|\{ v : (u,v) \in E_{HL} \}\right|^3
 &\geq \frac{1}{K^6 n^2}
     \left( \sum_{u : \beta \leq F(u) \leq K \beta}
              \left|\{ v : (u,v) \in E_{HL} \}\right| \right)^3
     \nonumber \\
 &\geq \frac{1}{K^6 n^2} \left|E_{HL}\right|^3
     \nonumber
\end{align}
and when $\sum_{(u,v) \in E_{HL}} F(u,v) \geq \frac12 \sum_{u,v} F(u,v)$
it follows from \eqref{eq:EHL} that this is at least
\begin{align}
  \frac{1}{8 K^9 n^2 \beta^3} \left(\sum_{u,v} F(u,v)\right)^3
&\gtrsim n \left( \mathbb E |S|^2\right)^{3/2},
    \nonumber
\end{align}
completing the proof.
\end{proof}

\section{Eigenvalues of graphs and surfaces}

\subsection{Graphs}
\label{sec:graphs}

We can now prove our main theorem.

\begin{theorem}
\label{thm:specificbounds}
If $G$ is an $n$-node graph, then for every $1 \leq k \leq n$, we have the following bounds.
  If $G$ is planar, then
  \begin{equation}\label{eq:planar}\lambda_k \leq O\left( d_{\max} \frac k n \right).\end{equation}
  If $G$ is of genus $g > 0$, then
  $$ \lambda_k \leq O\left( d_{\max}\frac k n   g (\log g)^2 \right). $$
  If $G$ is $K_h$-minor-free, then
  $$ \lambda_k \leq O\left( d_{\max}\frac k n h^6 \log h \right). $$
\end{theorem}

\begin{proof}
We prove the planar case; the other cases follow similarly.
Let $G=(V,E)$ be planar with maximum degree $d_{\max}$ and $n=|V|$.
First, by Theorem \ref{thm:eigenvalues}, we see that for any weight function $\omega : V \to \mathbb R_+$ and every $k \geq 1$, we have
\begin{equation}\label{eq:ff}
\lambda_{k} \lesssim \frac{d_{\max}}{\varepsilon^2 n} \beta^2,
\end{equation}
as long as $\frac{\beta^2}{\e^2} \leq \frac{n}{64}$,
where $\varepsilon = \varepsilon_{\lfloor n/4k\rfloor}(G,\omega)$
and $\beta=\beta_{\varepsilon/2}(V,\dist_{\omega})$.
Since $G$ is planar, by Theorem \ref{thm:KPR}, we have $\beta_{\varepsilon/2}(V,\dist_{\omega}) = O(1)$ for any $\omega$, hence
\eqref{eq:ff} implies
\begin{equation}\label{eq:primbound}
\lambda_{k} \lesssim \frac{d_{\max}}{\left(\varepsilon_{\lfloor n/4k\rfloor}(G,\omega)\right)^2\, n}.
\end{equation}
Using Corollaries \ref{cor:subset-crossings} and \ref{cor:ixn-basic}, we see that
for some constant $c_0 \geq 1$ and any $c_0 \leq r \leq |V|$, if $F \in \mathcal F_r(G)$, i.e. $F$ if a $\mu$-flow
with $\supp(\mu) \subseteq {V \choose r}$, then
$$
\vcon(F) \gtrsim \frac{r^5}{n}.
$$
Now, by Theorem \ref{thm:duality}, this implies that for $r \geq c_0$, there exists
a weight $\omega_r : V \to \mathbb R_+$ with
$\varepsilon_r(G,\omega_r) \gtrsim \frac{1}{r^2} \sqrt{r^5/n} = \sqrt{r/n}.$
In particular, for some constant $c_1 \geq c_0$, it follows that $r \geq c_1$
implies
\begin{equation}\label{eq:lastchange}
\frac{\beta^2}{\eps^2} \leq \frac{n}{64}\,.
\end{equation}

If $\lfloor n/8k \rfloor < c_1$, then \eqref{eq:planar} holds trivially using the bound $\lambda_k \leq 2\,d_{\max}$
for all $1 \leq k \leq n$.  Finally, if $r = \lfloor n/8k\rfloor \geq c_1$, then \eqref{eq:lastchange}
allows us to employ
\eqref{eq:primbound} to obtain,
$$
\lambda_{k} \lesssim \frac{d_{\max}}{\left(\varepsilon_{r}(G,\omega_r)\right)^2\, n} \lesssim \frac{d_{\max}}{r} \lesssim d_{\max}\,\frac{k}{n},
$$
completing the proof.

\begin{remark}[Asymptotic dependence on $k$]\label{rem:tightness}
We remark that the asymptotic dependence on $k$ in Theorem \ref{thm:specificbounds} is tight.
First, consider the eigenvalues $\lambda'_1 \leq \cdots \leq \lambda'_n$
for the $n$-node path graph $P_n$. It is a straightforward calculation to verify that the eigenvalues
are precisely the set $$\{ 2 - 2 \cos(2\pi k/n) : 1 \leq k \leq n/2 \},$$
and each such eigenvalue has multiplicity at most 2.
In particular, $\lambda'_k \asymp \frac{k^2}{n^2}$ for all $k \geq 2$.

Now, since the $n \times n$ grid graph $G_n$ is the Cartesian product graph $P_n \times P_n$,
it is easy to verify that the eigenvalues are precisely
$$
\{ \lambda_{i,j} = \lambda'_i + \lambda'_j : 1 \leq i,j \leq n \}.
$$
\end{remark}
In particular, since $\lambda_{i,j} \asymp \max(i^2,j^2)/n^2$,
we have $\lambda_k(G_n) \asymp \frac{k}{n^2} \asymp \frac{k}{|G_n|}.$
\end{proof}

Finally, we use the Robertson-Seymour structure theorem
to prove Theorem \ref{thm:minorclosed}.

\begin{proof}[Proof of Theorem \ref{thm:minorclosed}]
If $\mathcal F$ is any minor-closed family of graphs that
does not contain all graphs, then by the deep Robertson-Seymour
structure theory \cite{RS90}, there exists some number $h \in \mathbb N$
such that no graph in $\mathcal F$ has $K_h$ as a minor.  An application
of Theorem \ref{thm:specificbounds} finishes the proof.
\end{proof}

\subsection{Surfaces}\label{sec:surfaces}
\newcommand{\Lap}{\Delta}
\newcommand{\DLap}[1]{\Delta^{(#1)}}
\newcommand{\LM}{\Delta_M}
\newcommand{\DLM}[1]{\Delta^{(#1)}_M}
\newcommand{\bary}{\beta}
\newcommand{\dc}{d^c}
\newcommand{\deltac}{\delta^c}
\newcommand{\Lc}{\Delta^c}
\newcommand{\eign}{\lambda^{(n)}}
\newcommand{\dvol}{\,dV}
\newcommand{\area}{\mathrm{area}}
\newcommand{\Aeps}{A_{\epsilon}}

In this section, we shall show how our result implies a bound on the eigenvalues
of the Laplacian of a compact Riemannian surface.
\begin{theorem}\label{riemannbound}
Let $(M,g)$ be a compact, orientable Riemannian surface of genus $g$ and area
$A$, and let $\Lap_M$ be its Laplacian.  The $k^\text{th}$ smallest Neumann
eigenvalue of  $\Lap_M$ is at most $$O\left(k(g+1)\log^2 (g+1) / A\right).$$
\end{theorem}
Intuitively, this theorem follows by applying the eigenvalue bound for genus $g$ graphs from Theorem~\ref{thm:specificbounds} to a sequence of successively finer meshes that approximate $M$.

Our proof will begin with the combinatorial Hodge theory of Dodziuk~\cite{Dodziuk}, which produces a sequence of finite-dimensional operators $\DLM{1}, \DLM{2},\dots$ whose eigenvalues converge to those of $\LM$.
Unfortunately, the objects that this produces will not be the Laplacians of unweighted graphs of bounded degree.  However, we will show that, when applied to a sufficiently nice triangulation, the operators produced by Dodziuk's theory can be approximated well enough by such graph Laplacians to establish our desired result.

\subsubsection{The Whitney Map and Combinatorial Hodge Theory}

We begin by recalling the basic setup of Dodziuk's combinatorial Hodge theory~\cite{Dodziuk}.
Let $\chi: K \rightarrow M$ be a finite triangulation of $M$ with vertices
$p_1,\dots,p_n \in K$.
For all $q \in \mathbb N$,
let $L^2 \Lambda^q=L^2\Lambda(M)$ be the space of square
integrable $q$-forms on $M$, and let $C^q = C^q(K)$ be the space of real
simplicial cochains on $K$.
We will identify each simplex $\sigma$ of $K$ with the corresponding cochain,
which allows us to write elements of  $C^q(K)$ as formal sums of the
$q$-simplices in $K$. For any triangle
$\sigma \in K$,
we will use $\area(\sigma)$
and $\diam(\sigma)$  to denote that area and diameter of $\chi \sigma$ with
respect to the Riemannian metric on $M$.

For each $p_i$,
let
$\bary_i : K\rightarrow \mathbb{R}$
equal the $p_i^{\text{th}}$ barycentric coordinate on simplices in $\mathrm{St}(p_i)$, the open star of $p_i$,  and 0 on $K\setminus \mathrm{St}(p_i)$.  This lets us define barycentric coordinate functions $\mu_i = \chi^*\bary_i$ on $M$.

Let $\sigma = [p_{i_0},\dots, p_{i_q}]$ be a $q$-simplex in $K$ with $i_0 \leq \dots \leq i_q$.  We define the \emph{Whitney map} $W: C^q(K)\rightarrow L^2 \Lambda$ to
be the linear map that takes each such simplex to
$$W\sigma = q! \sum_{k=0}^q  (-1)^k \mu_{i_k} d\mu_{i_0} \wedge \dots \wedge  \widehat{d\mu_{i_k}} \wedge \dots \wedge d\mu_{i_q}.$$




Whitney~\cite{Whitney} showed that the above definition gives a well-defined element of
$L^2 \Lambda^q$, even though the $\mu_i$ are not differentiable on the
boundaries of top-dimensional simplices.


The Riemannian metric endows $L^2 \Lambda^q$ with the inner product
$$(f,g) = \int_M f \wedge *g, $$
where $*$ is the Hodge star operator.  Using the Whitney map, this lets us define an inner product on $C^q$ by setting
$$ (a,a')= (Wa, Wa')$$
for $a, a' \in C^q$.  Let $\dc$ be the simplicial coboundary operator. Dodziuk defined the combinatorial codifferential $\deltac$ to be the adjoint of $\dc$ with respect to this inner product, and he
defined the combinatorial Laplacian $\Lc_q: C^q \rightarrow C^q$ by
$$\Lc_q = \dc \deltac + \deltac \dc.$$
In the remainder of this paper, we will only use the Laplacian on functions,  which we will denote by $\Lc := \Lc_0$.

To obtain a sequence of successively finer triangulations, we will use
Whitney's \emph{standard subdivision} procedure~\cite{Whitney}.  For a complex $K$, this produces
a new complex $SK$ in which each $q$-dimensional simplex of $K$ is divided into $2^q$ smaller simplices.  In contrast to barycentric subdivision, it is constructed in a way that prevents the simplices from becoming arbitrarily poorly conditioned under repeated subdivision.


Let $S_0 K = K$, and inductively define $S_{n+1} K = S \left( S_n K\right)$.
Dodziuk showed the following convergence result about the discrete Laplacians on functions:%
\footnote{Dodziuk and Patodi~\cite{DodPat} later proved an analogous result for the Laplacians on $q$-forms, for arbitrary $q$.}
\begin{theorem}[Dodziuk]\label{Dodziuk}
Let $\eign_i$ be the $i^{\text{th}}$ smallest eigenvalue of $\Lc(S_n K)$, and let  $\lambda_i$ be the $i^{\text{th}}$ smallest eigenvalue of $\Lap_M$.
Then $\eign_i \rightarrow \lambda_i$ as $n \rightarrow \infty$.
\end{theorem}

\subsubsection{Relating the Combinatorial and Graph Laplacians}

To relate the combinatorial Laplacian to a graph Laplacian, we will construct a
triangulation in which all of the triangles
have approximately the same volume, are fairly flat, and have vertex angles
bounded away from 0.  We will then show that the eigenvalues of
combinatorial Laplacians arising from such a triangulation and its standard
subdivisions can be bounded in terms of those of the Laplacian of an unweighted graph of
bounded degree.

\begin{lemma}
\label{triangulation}
There exist strictly positive universal constants $C_1$, $C_2$, $C_3$, and
$\theta$ such that, for any $\epsilon>0$,
every compact Riemannian surface
$M$ has a triangulation $K$ with the following properties:
\begin{enumerate}
\item \label{property1} For every triangle $\sigma \in K$, $\diam (\sigma )< \epsilon$, the
interior angles of $\sigma$ all lie in  $[\theta ,\pi -\theta]$, and
\[
\frac{1}{C_{2}}\leq \frac{\area (\sigma)}{\diam(\sigma )^{2}}\leq C_{2}.
\]
\item \label{property2} For any two triangles $\sigma_{1},\sigma_{2}\in K$,
\[
1/C_1 \leq
\frac{\area(\sigma_{1})}{\area(\sigma_{2})}\leq C_{1},
\]
and
\[
1/C_1 \leq
\frac{\diam(\sigma_{1})}{\diam(\sigma_{2})}\leq C_{1}.
\]

\item \label{degbound} The edges of $K$ are embedded as geodesics, and every vertex of $K$ has degree at most $C_{3}$.
\end{enumerate}
Furthermore, these properties are satisfied by $S_{n}K$ for all $n\geq 0$.
\end{lemma}

\begin{proof}
The existence of such a triangulation is established by Buser, Sepp\"al\"a, and Silhol~\cite{BuserSS}, following an argument
originally due to Fejes T\'oth~\cite{Toth}.   They do not explicitly state the degree
bound, but it follows immediately from the fact that the angles are bounded away
from zero.  The fact that these properties remain true under
subdivision follows from the basic properties of standard subdivision given by
Whitney~\cite{Whitney}.
\end{proof}


\begin{proof}[Proof of Theorem~\ref{riemannbound}]
For a given $\epsilon$, let $K_{\epsilon}$ be a triangulation with the
properties guaranteed by Lemma~\ref{triangulation}, and let $G=(V,E)$ be the 1-skeleton of $K_{\epsilon}$.
Let $f:V\rightarrow\mathbb{R}$, and let $f_{i}=f(p_{i})$.  We will show that,
for sufficiently small $\epsilon$,
\begin{equation}\label{relateRayleigh}
\frac{(f,\Lc f)}{(f,f)} \lesssim \frac{|V|}{A}\frac{\|f\|_{\mathcal{L}_{G}}^{2}}{\|f\|_{2}^2}
\end{equation}
for all $f$, and that this remains true when $K_{\epsilon}$ is
replaced by
$S_{n}K_{\epsilon}$ for any $n$.  By the variational characterization of eigenvalues, this implies that
$
\lambda_k(\Lc)\lesssim \frac{|V|}{A}\lambda_{k}(\mathcal{L}_G)
$.
By applying Theorem~\ref{thm:specificbounds} to $\mathcal{L}_G$, we obtain
\[
\lambda_k (\Lc)\lesssim \frac{|V|}{A}\lambda_{k}(\mathcal{L}_G) \lesssim
\frac{|V|}{A}\frac{k (g+1) \log^{2} g}{|V|}
=  \frac{k (g+1) \log^{2} g}{A}.
\]
This bound remains true as we subdivide $K_{\epsilon}$, so Theorem~\ref{riemannbound} now follows from Theorem~\ref{Dodziuk}.
It thus suffices to prove equation~\eqref{relateRayleigh}.

Let $\sigma =[p_{i_{0}},p_{i_{1}},p_{i_{2}}]$ be a triangle in $K_{\epsilon}$.  We can
write the restriction of $Wf$ to $\sigma$ in barycentric coordinates as
\[
Wf|_{\sigma} = f_{1}\mu_{i_1}+f_2\mu_{i_2}+f_3\mu_{i_3}.
\]
%
When $\epsilon$ is sufficiently small
compared to the minimum  radius of curvature of $M$,
we have
\[
\int_{\sigma}\mu_{i}\mu_{j}\dvol = \begin{cases}
(1\pm o(1))\area(\sigma)/6 & \text{if $i=j$}\\
(1\pm o(1))\area(\sigma)/12 & \text{if $i\neq j$}
\end{cases},
\]
where $dV$ is the volume element on $M$, and the $o(1)$ indicates a function that goes to zero with $\epsilon$.

This gives
\begin{align*}
\int_{\sigma} (Wf) \wedge *(Wf)&=\int_{\sigma}(Wf)^{2}dV  \\
&=\int_{\sigma}(f_{1}\mu_{i_1}+f_2\mu_{i_2}+f_3\mu_{i_3})^{2}\dvol\\
&\asymp \frac{\area(\sigma)}{6}(f_1^2+f_2^2+f_3^2+f_{1} f_{2}+ f_{1}f_{3} + f_{2}f_{3})\\
&= \frac{\area(\sigma )}{12} \left(f_{1}^{2}+f_{2}^{2}+f_{3}^{2}+(f_{1}+f_{2}+f_{3})^{2}\right)\\
&\geq \frac{\area(\sigma )}{12} ( f_{1}^{2}+f_{2}^{2}+f_{3}^{2}).
\end{align*}

Let $\Aeps$ be the maximum area of a triangle in $K_{\epsilon}$.  Since  all
triangles  have the same area up to a multiplicative constant, and  each vertex
appears in only a constant number of triangles, summing this
over all of the triangles in $K_{\epsilon}$ gives
\begin{equation}\label{denom}
(f,f)=\int_{M}(Wf)\wedge*(Wf)\gtrsim \Aeps \sum_{i=1}^{n}f_{i}^{2}= \Aeps \| f
\|_{2}^{2}.
\end{equation}

When restricted to $\sigma$,  we have
\[
\dc f|_{\sigma}=(f_{1}-f_{0})[p_{i_{0}},p_{i_{1}}] +
(f_{2}-f_{1})[p_{i_{1}},p_{i_{2}}] + (f_{2}-f_{0})[p_{i_{0}},p_{i_{2}}],
\]
so
\begin{equation}\label{wdf}
W\dc f|_{\sigma} =  \left( f_{1}-f_{0} \right)\left(\mu_{i_0}d\mu_{i_1} -
\mu_{i_1}d\mu_{i_0}\right)
+  \left( f_{2}-f_{1} \right)\left(\mu_{i_1}d\mu_{i_2} -
\mu_{i_2}d\mu_{i_1}\right)
+  \left( f_{2}-f_{0} \right)\left(\mu_{i_0}d\mu_{i_2} -
\mu_{i_2}d\mu_{i_0}\right).
\end{equation}
By again assuming that $\epsilon$ is sufficiently small and using the fact that the triangles in $K_{\epsilon}$ are all well-conditioned, we obtain by a simple calculation the estimate
\[
\int_{\sigma }d\mu_{i_{k}} * d\mu_{i_{k}}\lesssim \left(\frac{1}{\diam(\sigma)}
\right)^{2}\cdot\area(\sigma)\asymp 1
\]
for each $k \in \{0,1,2\}$, where the asymptotic equality of the last two quantities follows from
property~\ref{property1} of Lemma~\ref{triangulation}.
Applying this and Cauchy-Schwartz to equation~(\ref{wdf}), and using the fact
that the $\mu_{i_{j}}$ are bounded above by 1, gives
\[
\int_{\sigma} (Wd^c f) \wedge *(Wd^c f) \lesssim
(f_{1}-f_{0})^{2}+(f_{2}-f_{1})^{2}+(f_{2}-f_{0})^{2}.
\]
Summing this over all of the triangles and using Lemma~\ref{triangulation} then yields
\begin{equation}\label{num}
(df,df) = \int_{M}(W\dc f)\wedge *(W\dc f) \lesssim  \sum_{(i,j)\in E}\left(f_{i}-f_{j} \right)^{2} =\|f\|^2_{\mathcal L}.
\end{equation}

The total area of $M$ equals $A$, and the area of each triangle is within a
constant factor of $\Aeps$, so
$|V|\asymp A/\Aeps $.  If we combine this with the inequalities
in~\eqref{denom} and~\eqref{num}, we obtain
$$
\frac{(f,\Lc f)}{(f,f)} =\frac{(df,df)}{(f,f)}
\lesssim
 \frac{ \|f\|_{\mathcal L}^{2}}
 {\Aeps \|  f \|_{2}^{2} }
 \asymp \frac{|V|}{A}\frac{ \|f\|_{\mathcal L}^{2}}
 {\|  f \|_{2}^{2} }. $$
This proves equation~\eqref{relateRayleigh}, which completes the proof of Theorem~\ref{riemannbound}.
\end{proof}

\bibliography{references}

\begin{thebibliography}{10}

\bibitem{acns}
M.~Ajtai, V.~Chv\'atal, M.~Newborn, and E.~Szemer\'edi.
\newblock Crossing-free subgraphs.
\newblock In A.~Kotzig, A.~Rosa, G.~Sabidussi, and J.~Turgeon, editors, {\em
  Theory and Practice of Combinatorics: A Collection of Articles Honoring Anton
  Kotzig on the Occasion of His Sixtieth Birthday.}, volume~12 of {\em Annals
  of discrete mathematics}. North-Holland, Amsterdam, 1982.

\bibitem{Alon}
N.~Alon.
\newblock Eigenvalues and expanders.
\newblock {\em Combinatorica}, 6(2):83--96, 1986.

\bibitem{AlonKahale}
N.~Alon and N.~Kahale.
\newblock A spectral technique for coloring random 3-colorable graphs.
\newblock {\em SIAM Journal on Computing}, 26:1733--1748, 1997.

\bibitem{AlonMilman}
N.~Alon and V.~D. Milman.
\newblock $\lambda_1$, isoperimetric inequalities for graphs, and
  superconcentrators.
\newblock {\em J. Comb. Theory Series B}, 38:73--88, 1985.

\bibitem{AlonSeymourThomas}
N.~Alon, P.~Seymour, and R.~Thomas.
\newblock A separator theorem for graphs with an excluded minor and its
  applications.
\newblock In {\em STOC '90: proceedings of the 22nd annual ACM Symposium on
  Theory of Computing}, pages 293--299. ACM, 1990.

\bibitem{AlpertYaoMoreBetter}
C.~J. Alpert and S.-Z. Yao.
\newblock Spectral partitioning: the more eigenvectors, the better.
\newblock In {\em DAC '95: Proceedings of the 32nd ACM/IEEE conference on
  Design automation}, pages 195--200. ACM, 1995.

\bibitem{AspvallGilbert}
B.~Aspvall and J.~R. Gilbert.
\newblock Graph coloring using eigenvalue decomposition.
\newblock Technical report, Ithaca, NY, USA, 1983.

\bibitem{BLR}
P.~Biswal, J.~R. Lee, and S.~Rao.
\newblock Eigenvalue bounds, spectral partitioning, and metrical deformations
  via flows.
\newblock {\em J. ACM}, 57(3), 2010.

\bibitem{convex-optimization}
S.~Boyd and L.~Vandenberghe.
\newblock {\em Convex optimization}.
\newblock Cambridge University Press, Cambridge, 2004.

\bibitem{BrinPage}
S.~Brin and L.~Page.
\newblock The anatomy of a large-scale hypertextual web search engine.
\newblock In {\em Proceedings of the seventh {International} {Wide} {Web}
  {Conference}}, 1998.

\bibitem{BuserSS}
P.~Buser, M.~Sepp\"al\"a, and R.~Silhol.
\newblock Triangulations and moduli spaces of {Riemann} surfaces with group
  actions.
\newblock {\em Manuscripta Mathematica}, 88(2):209--224, 1995.

\bibitem{ChanSchlagZien}
P.~K. Chan, M.~D.~F. Schlag, and J.~Y. Zien.
\newblock Spectral k-way ratio-cut partitioning and clustering.
\newblock In {\em DAC '93: Proceedings of the 30th international conference on
  Design automation}, pages 749--754. ACM, 1993.

\bibitem{ChanGilbertTeng}
T.~F. Chan, J.~R. Gilbert, and S.-H. Teng.
\newblock Geometric spectral partitioning.
\newblock Xerox PARC, Tech. Report, 1994.

\bibitem{Cheeger}
J.~Cheeger.
\newblock A lower bound for the smallest eigenvalue of the {L}aplacian.
\newblock In R.~C. Gunning, editor, {\em Problems in Analysis}, pages 195--199.
  Princeton University Press, 1970.

\bibitem{Chung}
F.~R. Chung.
\newblock Diameters and eigenvalues.
\newblock {\em J. Amer. Math. Soc.}, 2:187--196, 1989.

\bibitem{Dodziuk}
J.~Dodziuk.
\newblock Finite-difference approach to the {H}odge theory of harmonic forms.
\newblock {\em American Journal of Mathematics}, 98(1):79--104, 1976.

\bibitem{DodPat}
J.~Dodziuk and V.~K. Patodi.
\newblock Riemannian structures and triangulations of manifolds.
\newblock {\em Journal of the Indian Mathematical Society (N.S.}, 40:1--52,
  1976.

\bibitem{DH72}
W.~E. Donath and A.~J. Hoffman.
\newblock Algorithms for partitioning of graphs and computer logic based on
  eigenvectors of connection matrices.
\newblock {\em IBM Technical Disclosure Bulletin}, 15:938--944, 1972.

\bibitem{DH73}
W.~E. Donath and A.~J. Hoffman.
\newblock Lower bounds for the partitioning of graphs.
\newblock {\em J. Res. Develop.}, 17:420--425, 1973.

\bibitem{Spreading00}
G.~Even, J.~Naor, S.~Rao, and B.~Schieber.
\newblock Divide-and-conquer approximation algorithms via spreading metrics.
\newblock {\em J. ACM}, 47(4):585--616, 2000.

\bibitem{FT03}
J.~Fakcharoenphol and K.~Talwar.
\newblock An improved decomposition theorem for graphs excluding a fixed minor.
\newblock In {\em Proceedings of 6th Workshop on Approximation, Randomization,
  and Combinatorial Optimization}, volume 2764 of {\em Lecture Notes in
  Computer Science}, pages 36--46. Springer, 2003.

\bibitem{Fi73}
M.~Fiedler.
\newblock Algebraic connectivity of graphs.
\newblock {\em Czechoslovak Mathematical Journal}, 23(98):298--305, 1973.

\bibitem{Fi75a}
M.~Fiedler.
\newblock Eigenvectors of acyclic matrices.
\newblock {\em Czechoslovak Mathematical Journal}, 25(100):607--618, 1975.

\bibitem{Fi75b}
M.~Fiedler.
\newblock A property of eigenvectors of nonnegative symmetric matrices and its
  applications to graph theory.
\newblock {\em Czechoslovak Mathematical Journal}, 25(100):619--633, 1975.

\bibitem{GilbertHutchinsonTarjan}
J.~R. Gilbert, J.~P. Hutchinson, and R.~E. Tarjan.
\newblock A separation theorem for graphs of bounded genus.
\newblock {\em Journal of Algorithms}, 5:391--407, 1984.

\bibitem{GrigoryanYau}
A.~Grigor'yan and S.-T. Yau.
\newblock Decomposition of a metric space by capacitors.
\newblock In {\em Proceedings of Symposia in Pure Mathematics (Special volume
  dedicated to L. Nirenberg)}, volume~65, pages 39--75, 1999.

\bibitem{Hall}
K.~M. Hall.
\newblock An r-dimensional quadratic placement algorithm.
\newblock {\em Management Science}, 17:219--229, 1970.

\bibitem{Harary94}
F.~Harary.
\newblock {\em Graph Theory}.
\newblock Addison-Wesley, Reading, MA, 1994.

\bibitem{Hersch70}
J.~Hersch.
\newblock Quatre propri\'et\'es isop\'erim\'etriques de membranes sph\'eriques
  homog\`enes.
\newblock {\em C. R. Acad. Sci. Paris S\'er. A-B}, 270:A1645--A1648, 1970.

\bibitem{Kelner}
J.~A. Kelner.
\newblock Spectral partitioning, eigenvalue bounds, and circle packings for
  graphs of bounded genus.
\newblock {\em SIAM J. Comput.}, 35(4):882--902, 2006.

\bibitem{KPR}
P.~Klein, S.~A. Plotkin, and S.~Rao.
\newblock Excluded minors, network decomposition, and multicommodity flow.
\newblock In {\em STOC '93: Proceedings of the twenty-fifth annual ACM
  Symposium on Theory of Computing}, pages 682--690. ACM, 1993.

\bibitem{Kleinberg}
J.~M. Kleinberg.
\newblock Authoritative sources in a hyperlinked environment.
\newblock {\em Journal of the ACM}, 46:668--677, 1999.

\bibitem{Korevaar}
N.~Korevaar.
\newblock Upper bounds for eigenvalues of conformal metrics.
\newblock {\em J. Differential Geometry}, 37:73--93, 1993.

\bibitem{Kostochka}
A.~V. Kostochka.
\newblock The minimum hadwiger number for graphs with a given mean degree of
  vertices.
\newblock {\em Metody Diskret. Analiz.}, (38):37--58, 1982.

\bibitem{KLMN05}
R.~Krauthgamer, J.~R. Lee, M.~Mendel, and A.~Naor.
\newblock Measured descent: a new embedding method for finite metrics.
\newblock {\em Geom. Funct. Anal.}, 15(4):839--858, 2005.

\bibitem{LN05}
J.~R. Lee and A.~Naor.
\newblock Extending {L}ipschitz functions via random metric partitions.
\newblock {\em Invent. Math.}, 160(1):59--95, 2005.

\bibitem{LSgenus}
J.~R. Lee and A.~Sidiropoulos.
\newblock Genus and the geometry of the cut graph.
\newblock In {\em Proceedings of the 21st ACM-SIAM Symposium on Discrete
  Algorithms (SODA)}, pages 193--201, 2010.

\bibitem{Leighton}
F.~T. Leighton.
\newblock {\em Complexity Issues in VLSI}.
\newblock MIT Press, 1983.

\bibitem{LiptonTarjan}
R.~J. Lipton and R.~E. Tarjan.
\newblock A separator theorem for planar graphs.
\newblock {\em SIAM J. of Appl. Math.}, 36:177--189, 1979.

\bibitem{MillerTengThurstonVavasis}
G.~L. Miller, S.-H. Teng, W.~Thurston, and S.~A. Vavasis.
\newblock Finite element meshes and geometric separators.
\newblock {\em SIAM J. Scientific Computing}, 1996.

\bibitem{Rao}
S.~Rao.
\newblock Small distortion and volume preserving embeddings for planar and
  {E}uclidean metrics.
\newblock In {\em SCG '99: Proceedings of the fifteenth annual Symposium on
  Computational Geometry}, pages 300--306. ACM, 1999.

\bibitem{RS90}
N.~Robertson and P.~D. Seymour.
\newblock Graph minors. {VIII}.\ {A} {K}uratowski theorem for general surfaces.
\newblock {\em J. Combin. Theory Ser. B}, 48(2):255--288, 1990.

\bibitem{RS04}
N.~Robertson and P.~D. Seymour.
\newblock Graph minors. {XX}. {W}agner's conjecture.
\newblock {\em J. Combin. Theory Ser. B}, 92(2):325--357, 2004.

\bibitem{ShiMalik}
J.~Shi and J.~Malik.
\newblock Normalized cuts and image segmentation.
\newblock {\em IEEE Trans. Pattern Anal. Mach. Intell.}, 22(8):888--905, 2000.

\bibitem{SinclairJerrum}
A.~J. Sinclair and M.~R. Jerrum.
\newblock Approximative counting, uniform generation and rapidly mixing
  {Markov} chains.
\newblock {\em Information and Computation}, 82(1):93--133, 1989.

\bibitem{SpielmanTeng}
D.~A. Spielman and S.-H. Teng.
\newblock Spectral partitioning works: planar graphs and finite element meshes.
\newblock {\em Linear Algebra Appl.}, 421(2-3):284--305, 2007.

\bibitem{Thomason}
A.~Thomason.
\newblock An extremal function for contractions of graphs.
\newblock {\em Math. Proc. Cambridge Philos. Soc.}, 95(2):261--265, 1984.

\bibitem{TolliverMiller}
D.~A. Tolliver and G.~L. Miller.
\newblock Graph partitioning by spectral rounding: Applications in image
  segmentation and clustering.
\newblock In {\em CVPR '06: Proceedings of the 2006 IEEE Computer Society
  Conference on Computer Vision and Pattern Recognition}, pages 1053--1060.
  IEEE Computer Society, 2006.

\bibitem{Toth}
L.~F. T\'oth.
\newblock Kreisausfüllungen der hyperbolischen ebene.
\newblock {\em Acta Mathematica Hungarica}, 4(1):103--110, 1953.

\bibitem{Wagner}
K.~Wagner.
\newblock \"{U}ber eine {E}igenschaft der ebenen {K}omplexe.
\newblock {\em Math. Ann.}, 114(1):570--590, 1937.

\bibitem{Whitney}
H.~Whitney.
\newblock {\em Geometric integration theory}.
\newblock Princeton University Press, 1957.

\bibitem{YangYau}
P.~Yang and S.~T. Yau.
\newblock Eigenvalues of the {Laplacian} of compact {Riemann} surfaces and
  minimal submanifolds.
\newblock {\em Ann. Scuola Norm. Sup. Pisa Cl. Sci.}, 4(7):55--63, 1980.

\end{thebibliography}
\bibliographystyle{abbrv}

\remove{
\appendix

\section{Additional proofs}

\subsection{The duality proof}

Here we present the proof of Theorem~\ref{thm:duality}.

\parskip1em \noindent
{\bf Theorem \ref{thm:duality} {\rm (Duality)}.}  {\it
  Let $G = (V, E)$ be a graph and let $r \leq |V|$.
  Then
  $$
  \max \left\{ \vphantom{\bigoplus} \varepsilon_r(G,\omega)
              \middle| \omega : V \to \mathbb R_+ \right\}
  = \frac{1}{r^2} \min \left\{ \vphantom{\bigoplus} \sqrt{\vcon(F)}
              \middle| F \in \mathcal F_r(G) \right\}
  $$
}

\subsection{Congestion measures}
\label{sec:measures}

\subsection{Conclusion of the heavy edges, light endpoints analysis}
\label{sec:heavy}

}

\end{document}